\documentclass[10pt]{article}

\usepackage{latexsym}
\usepackage{amsfonts}
\usepackage{graphicx}
\usepackage{amsmath}
\usepackage{amssymb}
\usepackage{geometry}
\usepackage[latin1]{inputenc}
\usepackage{url}        %ajout
\usepackage[thmmarks,amsmath]{ntheorem}  %ajout
\usepackage{subfigure}
\usepackage[colorlinks=true]{hyperref}

\makeatletter \@addtoreset{equation}{section}

\makeatletter

\def\ph{\varphi}

\def\be{\begin{equation}}   \def\ee{\end{equation}}
\def\ba   {\begin{array}}      \def\ea   {\end{array}}
\def\bea  {\begin{eqnarray}}   \def\eea  {\end{eqnarray}}
\def\bean {\begin{eqnarray*}}  \def\eean {\end{eqnarray*}}

\newtheorem{theorem} {Theorem}
\newtheorem{lemma}{Lemma}
\newtheorem{definition} {Definition}

\newtheorem{corollary} {Corollary}
\newtheorem{remark}{Remark}

\newtheorem{conjecture} {Conjecture}

\newcommand{\pre}{\mathrm{Re}}
\newcommand{\pim}{\mathrm{Im}}
\newcommand{\bi} {\ensuremath{{\bf i}}}
\newcommand{\bo} {\ensuremath{{\bf i_1}}}
\newcommand{\bos}{\ensuremath{{\bf i_1^{\text 2}}}}
\newcommand{\bts}{\ensuremath{{\bf i_2^{\text 2}}}}

\newcommand{\bj}{\ensuremath{{\bf j}}}

\newcommand {\bjp}{\ensuremath{{\bf j_1}}}
\newcommand {\bjps}{\ensuremath{{\bf j_1^{\text 2}}}}
\newcommand {\bjd}{\ensuremath{{\bf j_2}}}
\newcommand {\bjds}{\ensuremath{{\bf j_2^{\text 2}}}}

\newcommand{\bt} {\ensuremath{{\bf i_2}}}
\newcommand{\bb} {\ensuremath{{\bf i_3}}}
\newcommand{\bbs} {\ensuremath{{\bf i_3^{\text 2}}}}
\newcommand{\bq} {\ensuremath{{\bf i_4}}}
\newcommand{\bqs} {\ensuremath{{\bf i_4^{\text 2}}}}
\newcommand{\bn} {\ensuremath{{\bf i_n}}}

\newcommand{\bnm} {\ensuremath{{\bf i_{n-1}}}}

\newcommand{\bk} {\ensuremath{{\bf i_{k}}}}
\newcommand {\bjt}{\ensuremath{{\bf j_3}}}
\newcommand {\bjts}{\ensuremath{{\bf j_3^{\text 2}}}}
\newcommand {\bjk}{\ensuremath{{\bf j_k}}} %AJOUT
\newcommand{\bil} {\ensuremath{{\bf i_l}}}
\newcommand{\bim} {\ensuremath{{\bf i_m}}}
\newcommand{\eb} {\ensuremath{\gamma_1}}

\newcommand{\ett} {\ensuremath{\gamma_2}}
\newcommand{\etc} {\ensuremath{\Ol{\gamma}_2}}

\newcommand{\Ol}{\overline}

\newcommand{\mC}{\ensuremath{\mathbb{C}}}
\newcommand{\mD}{\ensuremath{\mathbb{D}}}
\newcommand{\mN}{\ensuremath{\mathbb{N}}}

\newcommand{\mR}{\ensuremath{\mathbb{R}}}
\newcommand{\mT}{\ensuremath{\mathbb{T}}}

\newcommand{\mH}{\ensuremath{\mathbb{D}}}

\newcommand{\mM}{\ensuremath{\mathbb{M}}}
\newcommand{\mMan}{\ensuremath{\mathcal{M}}}
\newcommand{\mTet}{\ensuremath{\mathcal{T}}}

\newcommand{\oa}{\left\{}
\newcommand{\fa}{\right\}}

\newcommand{\Qpit}{\ensuremath{\left\lbrace Q_{p,c}^m(0) \right\rbrace_{m=1}^{\infty}}}
\newcommand{\Qit}{\ensuremath{\left\lbrace Q_{3,c}^m(0) \right\rbrace_{m=1}^{\infty}}}
\newcommand{\mManp}{\ensuremath{\mathcal{M}^p}}
\newcommand{\mManth}{\ensuremath{\mathcal{M}^3}}

\newcommand{\mHt}{\ensuremath{\mathcal{H}^2}}
\newcommand{\mHth}{\ensuremath{\mathcal{H}^3}}
\newcommand{\mHyb}{\ensuremath{\mathcal{H}}}

\newcommand{\cc}{\ensuremath{\overline{c}}}
\newcommand{\bHpc}{\ensuremath{\mathbf{H}_{p,c}}}
\newcommand{\vecxy}{\ensuremath{\begin{pmatrix}
x \\ y
\end{pmatrix}}}

\newcommand{\vecab}{\ensuremath{\begin{pmatrix}
a \\ b
\end{pmatrix}}}

{\theoremsymbol{$\square$}
\theorembodyfont{}
\newtheorem*{demo}{\textbf{Proof.}}}

\newcommand{\bjl}{\ensuremath{\bf j_l}}
\newcommand{\bjm}{\ensuremath{\bf j_m}}
\newcommand{\biq}{\ensuremath{\bf i_q}}
\newcommand{\bis}{\ensuremath{\bf i_s}}

\begin{document}

%Infos article et auteurs première page
\markboth{P.-O. Parisé and D. Rochon}{A Study of Dynamics of the Tricomplex Polynomial $\eta^p+c$}

\title{A Study of Dynamics of the Tricomplex Polynomial $\eta^p+c$
}
\author{Pierre-Olivier Parisé\thanks{E-mail: \texttt{Pierre-Olivier.Parise@uqtr.ca}} \and
Dominic Rochon\thanks{E-mail: \texttt{Dominic.Rochon@uqtr.ca}}}
%\date{\today}
\date{Département de mathématiques et
d'informatique \\ Université du Québec à Trois-Rivières \\
C.P. 500 Trois-Rivières, Québec \\ Canada, G9A 5H7}

\maketitle

\vspace{1cm} \noindent\textbf{AMS subject classification:} 37F50, 32A30, 30G35, 00A69

%-------------------------------------------------------------------------------

\begin{abstract}
In this article, we give the exact interval of the cross section of the so-called \textit{Mandelbric} set generated by the polynomial $z^3+c$ where $z$ and $c$ are complex numbers. Following that result, we show that the \textit{Mandelbric} defined on the hyperbolic numbers $\mH$ is a square with its center at the origin. Moreover, we define the \textit{Multibrot} sets generated by a polynomial of the form $Q_{p,c}(\eta )=\eta^p+c$ ($p \in \mN$ and $p \geq 2$) for tricomplex numbers. More precisely, we prove that the tricomplex \textit{Mandelbric} has four principal slices instead of eight principal 3D slices that arise for the case of the tricomplex Mandelbrot set. Finally, we prove that one of these four slices is an octahedron.
\end{abstract}

\section{Introduction}
In 1982, A. Douady and J. H. Hubbard \cite{Dou} studied dynamical systems generated by iterations of the quadratic polynomial $z^2+c$. One of the main results of their work was the proof that the well-known Mandelbrot set for complex numbers is a connected set. It is also well known that the Mandelbrot set is bounded by a discus of radius 2 and crosses the real axis on the interval $\left[ -2 , \frac{1}{4} \right]$. Considering functions of the form $z^p+c$ where $z,c$ are complex numbers and $p$ may be an integer, and a rational or a real number being greater than 2, T. V. Papathomas, B. Julesz, U. G. Gujar and  V. G. Bhavsar \cite{Papathomas, Gujar} explored the sets generated by these functions called \textit{Multibrot} sets. The last two authors remark that these \textit{polynomials} generate rich fractal structures. This was the starting point for other researchers such as D. Schleicher (\cite{Di,Lau}), E. Lau \cite{Lau}, X. Sheng and M. J. Spurr \cite{Sheng}, X.-D. Liu et al. \cite{chine6} and many others to study symmetries in \textit{Multibrot} sets, their connectivity and the discus that bound these sets.

In 1990, P. Senn \cite{Senn} suggested to define the Mandelbrot set for another set of numbers: the hyperbolic numbers (also called duplex numbers). He remarked that the Mandelbrot set for this number structure seemed to be a square. Four years later, a proof of this statement was given by W. Metzler \cite{MET}.

Unless these works were done in the complex plane, so in 2D, several mathematicians questioned themselves on a generalization of the Mandelbrot set in three dimensions.
In 1982, A. Norton \cite{NORTON} succeeded to bring a first method to visualize fractals in 3D using the quaternions. In 2000, D. Rochon \cite{Rochon1} used the bicomplex numbers set $\mM (2)$ to give a 4D definition of the so-called Mandelbrot set and made 3D slices to get the \textit{Tetrabrot}. Later, X.-y. Wang and W.-j. Song \cite{Chine1} follow D. Rochon's work to establish the \textit{Multibrots} sets for bicomplex numbers. Several years before, D. Rochon and V. Garant-Pelletier (\cite{GarantRochon}, \cite{GarantPelletier}) gave a definition of the Mandelbrot set for multicomplex numbers denoted by $\mMan_n$ and restrict their explorations to the tricomplex case. Particularly, they found eight principal slices of the tricomplex Mandelbrot set and proved that one of these 3D slices, typically named the \textit{Perplexbrot}, is an octahedron of edges $\frac{9\sqrt{2}}{8}$.

In this article, we investigate the \textit{Multibrot} sets for complex, hyperbolic and tricomplex numbers, respectively, denoted by $\mMan^p$, $\mHyb^p$ and $\mMan_3^p$ when $p$ is an integer greater than one. We emphasize on the sets $\mMan^3$, $\mHyb^3$ and $\mMan_3^3$ respectively called the \textit{Mandelbric}, \textit{Hyperbric} and tricomplex \textit{Mandelbric}. Precisely, in the second section, we recall some definitions and properties of tricomplex numbers denoted by $\mM (3)$. We remark that complex and hyperbolic numbers are embedded in $\mM (3)$ and also that there are other interesting subsets of $\mM (3)$. In the third section, we review the definition and properties of \textit{Multibrots} sets. Particularly, we prove that the set $\mMan^3$ crosses the real axis on $\left[- \frac{2}{3\sqrt{3}}, \frac{2}{3\sqrt{3}} \right]$. In section four, we define \textit{Multibrot} sets for hyperbolic numbers. Particularly, based on Metzler's article (see \cite{MET}), we prove that $\mHyb^3$ is a square where its center is the origin. Finally, in the fifth section, we define the tricomplex \textit{Multibrot} sets corresponding to the polynomial $\eta^p+c$ where $\eta$ and $c$ are tricomplex numbers and $p\geq 2$ is an integer. We obtain, for the case where $p=3$, that there are four principal 3D slices of $\mMan_3^3$ instead of eight like it is showed in \cite{GarantRochon} for $\mMan_3^2$. Moreover, we prove that one of these four slices, typically named the \textit{Perplexbric}, is an octahedron of edges $\frac{2\sqrt{2}}{3\sqrt{3}}$.

\section{Tricomplex numbers}
In this section, we begin by a short introduction of the tricomplex numbers space $\mM (3)$. One may be refer to \cite{Baley}, \cite{GarantPelletier} and \cite{Vajiac} for more details on the next properties.

A tricomplex number $\eta$ is composed of two coupled bicomplex numbers $\zeta_1$, $\zeta_2$ and an imaginary unit $\bb$ such that
\begin{equation}
\eta=\zeta_1 + \zeta_2 \bb \label{eq2.1}
\end{equation}
where $\bbs=-1$. The set of such tricomplex numbers is denoted by $\mM (3)$. Since $\zeta_1,\zeta_2 \in \mM (2)$, we can write them as $\zeta_1=z_1+ z_2\bt$ and $\zeta_2=z_3+ z_4\bt$ where $z_1,z_2,z_3,z_4 \in \mM (1)\simeq \mC$. In that way, \eqref{eq2.1} can be rewritten as
\begin{equation}
\eta=z_1+ z_2 \bt+  z_3 \bb+  z_4 \bjt\label{eq2.2}
\end{equation}
where $\bts=-1$, $\bt \bb = \bb \bt = \bjt$ and $\bjts=1$. Moreover, as $z_1$, $z_2$, $z_3$ and $z_4$ are complex numbers, we can write the number $\eta$ in a third form as
\begin{align}
\eta&=a+ b\bo + (c+  d\bo)\bt + (e +  f\bo)\bb + (g +  h\bo)\bjt\notag\\
&=a+ b\bo +  c\bt +  d\bjp +  e\bb +  f\bjd + g \bjt +  h\bq\label{eq2.3}
\end{align}
where $\bos=\bqs=-1$, $\bq =\bo \bjt = \bo \bt \bb$, $\bjd = \bo \bb = \bb \bo$, $\bjds=1$, $\bjp=\bo \bt = \bt \bo$ and $\bjps=1$. After ordering each term of \eqref{eq2.3}, we get the following representations of the set of tricomplex numbers:
\begin{align}
\mM (3) &:= \oa \eta = \zeta_1 +  \zeta_2\bb \, |\, \zeta_1, \zeta_2 \in \mM (2) \fa \notag\\
&=\oa z_1+ z_2 \bt+  z_3 \bb+ z_4 \bjt  \, |\, z_1,z_2,z_3,z_4 \in \mM (1) \fa \notag\\
&=\oa x_0+ x_1\bo + x_2\bt  +  x_3\bb + x_4\bq  +  x_5\bjp +  x_6 \bjd+  x_7\bjt \, |\, x_i \in \mM (0)=\mR \right. \\
& \left. \qquad \qquad \text{ for } i=0,1,2, \ldots , 7 \fa \text{.} \label{EqRep}
\end{align}
Let $\eta_1=\zeta_1+\zeta_2\bb$ and $\eta_2=\zeta_3 + \zeta_4\bb$ be two tricomplex numbers with $\zeta_1,\zeta_2,\zeta_3,\zeta_4 \in \mM (2)$. We define the equality, the addition and the multiplication of two tricomplex numbers as
\begin{align}
\eta_1&=\eta_2 \text{ iff } \zeta_1=\zeta_3 \text{ and } \zeta_2=\zeta_4 \label{eq2.4}\\
\eta_1 + \eta_2 &:= (\zeta_1 + \zeta_3) + (\zeta_2+\zeta_4)\bb \label{eq2.5}\\
\eta_1 \cdot \eta_2&:= (\zeta_1\zeta_3-\zeta_2\zeta_4)+(\zeta_1\zeta_4 + \zeta_2\zeta_3)\bb \label{eq2.6}\text{.}
\end{align}
Table \ref{tabC1} shows the results after multiplying each tricomplex imaginary unity two by two.
\begin{table}
\centering
\begin{tabular}{c|*{9}{c}}
$\cdot$ & 1 & $\mathbf{i_1}$ & $\mathbf{i_2}$ & $\mathbf{i_3}$ & $\mathbf{i_4}$ & $\mathbf{j_1}$ & $\mathbf{j_2}$ & $\mathbf{j_3}$\\\hline
1 & 1 & $\mathbf{i_1}$ & $\mathbf{i_2}$ & $\mathbf{i_3}$ & $\mathbf{i_4}$ & $\mathbf{j_1}$ & $\mathbf{j_2}$ & $\mathbf{j_3}$\\
$\mathbf{i_1}$ & $\mathbf{i_1}$ & $-\mathbf{1}$ & $\mathbf{j_1}$ & $\mathbf{j_2}$ & $-\mathbf{j_3}$ & $-\mathbf{i_2}$ & $-\mathbf{i_3}$ & $\mathbf{i_4}$\\
$\mathbf{i_2}$ & $\mathbf{i_2}$ & $\mathbf{j_1}$ & $-\mathbf{1}$ & $\mathbf{j_3}$ & $-\mathbf{j_2}$ & $-\mathbf{i_1}$ & $\mathbf{i_4}$ & $-\mathbf{i_3}$\\
$\mathbf{i_3}$ & $\mathbf{i_3}$ & $\mathbf{j_2}$ & $\mathbf{j_3}$ & $-\mathbf{1}$ & $-\mathbf{j_1}$ & $\mathbf{i_4}$ & $-\mathbf{i_1}$ & $-\mathbf{i_2}$\\
$\mathbf{i_4}$ & $\mathbf{i_4}$ & $-\mathbf{j_3}$  & $-\mathbf{j_2}$ & $-\mathbf{j_1}$ & $-\mathbf{1}$ & $\mathbf{i_3}$ & $\mathbf{i_2}$ & $\mathbf{i_1}$\\
$\mathbf{j_1}$ & $\mathbf{j_1}$ & $-\mathbf{i_2}$  & $-\mathbf{i_1}$ & $\mathbf{i_4}$ & $\mathbf{i_3}$ & $\mathbf{1}$ & $-\mathbf{j_3}$ & $-\mathbf{j_2}$\\
$\mathbf{j_2}$ & $\mathbf{j_2}$ & $-\mathbf{i_3}$  & $\mathbf{i_4}$ & $-\mathbf{i_1}$ & $\mathbf{i_2}$ & $-\mathbf{j_3}$ &  $\mathbf{1}$ & $-\mathbf{j_1}$\\
$\mathbf{j_3}$ & $\mathbf{j_3}$ & $\mathbf{i_4}$ &$-\mathbf{i_3}$  & $-\mathbf{i_2}$ & $\mathbf{i_1}$ & $-\mathbf{j_2}$ & $-\mathbf{j_1}$ &  $\mathbf{1}$ \\
\end{tabular}
\caption{Products  of tricomplex imaginary units}\label{tabC1}
\end{table}
The set of tricomplex numbers with addition $+$ and multiplication $\cdot$ forms a commutative ring with zero divisors.

A tricomplex number has a useful representation using the idempotent elements $\ett =\frac{1+\bjt}{2}$ and $\etc =\frac{1-\bjt}{2}$. Recalling that $\eta = \zeta_1 +  \zeta_2\bb$ with $\zeta_1, \zeta_2 \in \mM (2)$, the idempotent representation of $\eta$ is
\begin{equation}
\eta = (\zeta_1- \zeta_2\bt)\ett + (\zeta_1+ \zeta_2\bt)\etc \label{eq2.7}\text{.}
\end{equation}
The representation \eqref{eq2.7} of a tricomplex number is useful to add and multiply tricomplex numbers because it allows to do these operations term-by-term. In fact, we have the following theorem (see \cite{Baley}):
\begin{theorem}\label{theo2.2}
Let $\eta_1=\zeta_1 +  \zeta_2\bb$ and $\eta_2=\zeta_3 +  \zeta_4\bb$ be two tricomplex numbers. Let $\eta_1=u_1\ett + u_2 \etc$ and $\eta_2=u_3\ett + u_4\etc$ be the idempotent representation \eqref{eq2.7} of $\eta_1$ and $\eta_2$. Then,
\begin{enumerate}
\item $\eta_1+\eta_2=(u_1+u_3)\ett + (u_2+u_4)\etc$;
\item $\eta_1 \cdot \eta_2 = (u_1 \cdot u_3)\ett + (u_2 \cdot u_4)\etc$;
\item $\eta_1^m=u_1^m \ett + u_2^m \etc$ $\forall m \in \mN$.
\end{enumerate}
\end{theorem}
Moreover, we define a $\mM (3)$-\textit{cartesian} set $X$ of two subsets $X_1,X_2\subseteq\mM (2)$ as follows:
\begin{align}
X=X_1\times_{\ett}X_2:=\oa \eta =\zeta_1+ \zeta_2\bb \in \mM (3) \, | \, \eta = u_1 \ett + u_2 \etc , u_1 \in X_1 \text{ and } u_2 \in X_2 \fa \text{.}\label{eq2.14}
\end{align}

Let define the norm $\Vert \cdot \Vert_3 :\, \mM (3) \rightarrow \mR$ of a tricomplex number $\eta=\zeta_1 + \zeta_2\bb$ as
\begin{align}
\Vert \eta \Vert_3 & := \sqrt{\Vert\zeta_1\Vert_2^2+\Vert\zeta_2\Vert_2^2}=\sqrt{\sum_{i=1}^2|z_i|^2+\sum_{i=3}^4|z_i|^2}\label{eq2.15}\\
&=\sqrt{\sum_{i=0}^7x_i^2}.\notag
\end{align}
According to the norm \eqref{eq2.15}, we say that a sequence $\oa s_m \fa_{m=1}^{\infty} $ of tricomplex numbers is bounded if and only if there exists a real number $M$ such that $\Vert s_m \Vert_3 \leq M$ for all $m \in \mN$. Now,
according to \eqref{eq2.14}, we define two kinds of tricomplex discus:
\begin{definition}\label{def2.1}
Let $\alpha = \alpha_1+\alpha_2 \bb \in \mM (3)$ and set $\bf{r_2}\geq \bf{r_1} >0$.
\begin{enumerate}
\item The open discus is the set
\begin{align}
\bf{D_3}(\alpha ; \bf{r_1},\bf{r_2})&:= \oa  \eta \in \mM (3) \, | \, \eta =\zeta_1 \ett + \zeta_2\etc , \, \Vert \zeta_1-(\alpha_1- \alpha_2 \bt)\Vert_2 <\bf{r_1} \text{ and }\right.\notag \\ &\left. \qquad \qquad \qquad \Vert \zeta_2 - (\alpha_1+\alpha_2 \bt) \Vert_2 <\bf{r_2} \fa \text{.}\label{eq2.141}
\end{align}
\item The closed discus is the set
\begin{align}
\Ol{\bf{D_3}}(\alpha ; \bf{r_1},\bf{r_2})&:= \oa  \eta \in \mM (3) \, | \, \eta =\zeta_1 \ett + \zeta_2\etc , \, \Vert \zeta_1-(\alpha_1- \alpha_2\bt)\Vert_2 \leq \bf{r_1} \text{ and }\right.\notag \\ &\left. \qquad \qquad \qquad \Vert \zeta_2 - (\alpha_1+ \alpha_2 \bt) \Vert_2 \leq \bf{r_2} \fa \text{.}\label{eq2.142}
\end{align}
\end{enumerate}
\end{definition}

We end this section by several remarks about subsets of $\mM (3)$. Let the set $\mC (\bk ):=\oa \eta = x_0 + x_1 \bk \, | \, x_0, x_1 \in \mR \fa, \bk \in \oa \bo , \bt , \bb , \bq \fa $. So, $\mC (\bk )$ is a subset of $\mM (3)$ for $k=1,2,3,4$ and we also remark that they are all isomorphic to $\mC$. Furthermore, the set 
$$\mH (\bjk ):=\oa x_0 + x_1 \bjk \, | \, x_0, x_1 \in \mR\fa$$ where  $\bjk \in \oa \bjp , \bjd, \bjt  \fa$ is a subset of $\mM (3)$ and is isomorphic to the set of hyperbolic numbers $\mH$ for $k\in\{1,2,3,4\}$ (see \cite{RochonShapiro, vajiac2} and \cite{Sobczyk} for further details about the set $\mH$ of hyperbolic numbers). Moreover, we define three particular subsets of $\mM (3)$ (see \cite{GarantRochon} and \cite{Parise}).
\begin{definition} \label{Mikil}
Let $\bk , \bil \in \oa 1 , \bo , \bt , \bb , \bq , \bjp , \bjd , \bjt \fa$ where $\bk \neq \bil$. The first set is
\begin{equation}
\mM (\bk , \bil ):= \oa x_1 + x_2\bk + x_3\bil + x_4\bk\bil \, | \, x_i \in \mR, i=1, \ldots, 4 \fa \text{.}
\end{equation}
\end{definition}
It is easy to see that $\mM (\bk , \bil )$ is closed under addition and multiplication of tricomplex numbers and that $\mM (\bk , \bil ) \simeq \mM (2)$ except for
the \textit{biduplex} sets $\mM (\bjp,\bjd )$, $\mM (\bjp , \bjt )$ and $\mM (\bjd , \bjt )$ (see \cite{GarantRochon}). 
\begin{definition}\label{Mikilim}
Let $\bk , \bil , \bim \in \oa \bo , \bt , \bb , \bq , \bjp , \bjd , \bjt \fa$ with $\bk \neq \bil$, $\bk \neq \bim$ and $\bil \neq \bim$. The second subset is
\begin{equation}
\mM (\bk , \bil , \bim ):= \oa x_1\bk + x_2\bil + x_3\bim + x_4\bk\bil\bim \, | \, x_i \in \mR, i=1, \ldots, 4 \fa \text{.}
\end{equation}
\end{definition}
Using Table \ref{tabC1}, we can easily verify that for any tricomplex number $\zeta \in \mM (\bk , \bil , \bim )$, $\zeta^3 \in \mM (\bk , \bil , \bim )$. In section \ref{sec5}, this fact will be useful to characterize some principal 3D slices of the tricomplex \textit{Mandelbric}.
\begin{definition}\label{Tikilim}
Let $\bk , \bil , \bim \in \oa 1, \bo , \bt , \bb , \bq , \bjp , \bjd , \bjt \fa$ with $\bk \neq \bil$, $\bk \neq \bim$ and $\bil \neq \bim$. The third subset is
\begin{equation}
\mT (\bim , \bk , \bil ):= \oa x_1\bk + x_2\bil + x_3\bim \, | \, x_1,x_2,x_3 \in \mR \fa \text{.}
\end{equation}
\end{definition}
This last set is not closed under multiplication depending on the number of times $k$ you multiply a tricomplex number in this set. For $k$ even, two cases may occur depending on the choice of the tricomplex imaginary units: the case that $\mT (\bim , \bk , \bil ) \subseteq \mM (\bk , \bil )$ or the case that the result of multiplying tricomplex numbers in $\mT (\bim , \bk , \bil )$ is not closed in the set $\mM (\bk , \bil )$. The first case arises if one of the imaginary unit $\bk , \bil , \bim$ is 1 or if the product $\bk\bil = \pm \bim$. Whenever one of these conditions are not fulfilled, the result is not closed in the set $\mM (\bk , \bil )$. On the other hand, if $k$ is odd, the first case stay the same but the second is always closed in the set $\mM (\bk , \bil , \bim )$. These facts are direct consequences of the definition of the tricomplex imaginary units.

\section{Generalized Mandelbrot sets for complex numbers}
In this section, we investigate \textit{Multibrot} sets and recall some of their properties that come from \cite{Gujar,chine6,2Noah,Parise,Chine1}. We also obtain some specific results for the \textit{Mandelbric} set $\mMan^3$ generated by the complex polynomial $Q_{3,c}(z)=z^3+c$.

\subsection{\textit{Multibrot} sets}
\begin{figure}
\subfigure[$\mathcal{M}^3$ set: $-1.5\leq \pre (c) \leq 1.5$ and $-1.5 \leq \pim (c) \leq 1.5$]{%
\includegraphics[scale=0.225]{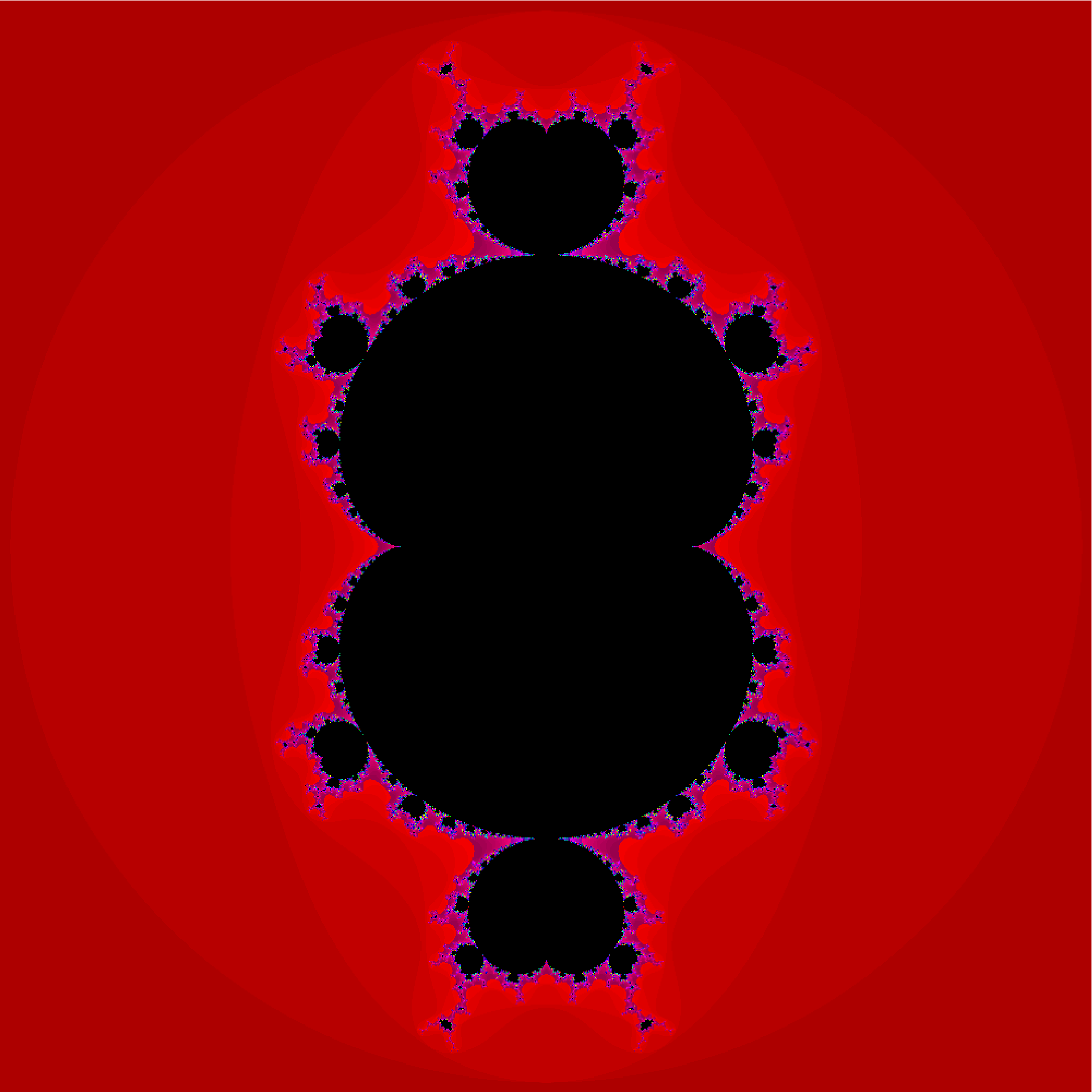}\label{fig2.1}}
\subfigure[$\mathcal{M}^4$ set: $-1.3\leq \pre (c) \leq 1.3$ and $-1.3 \leq \pim (c) \leq 1.3$]{%
\includegraphics[scale=0.225]{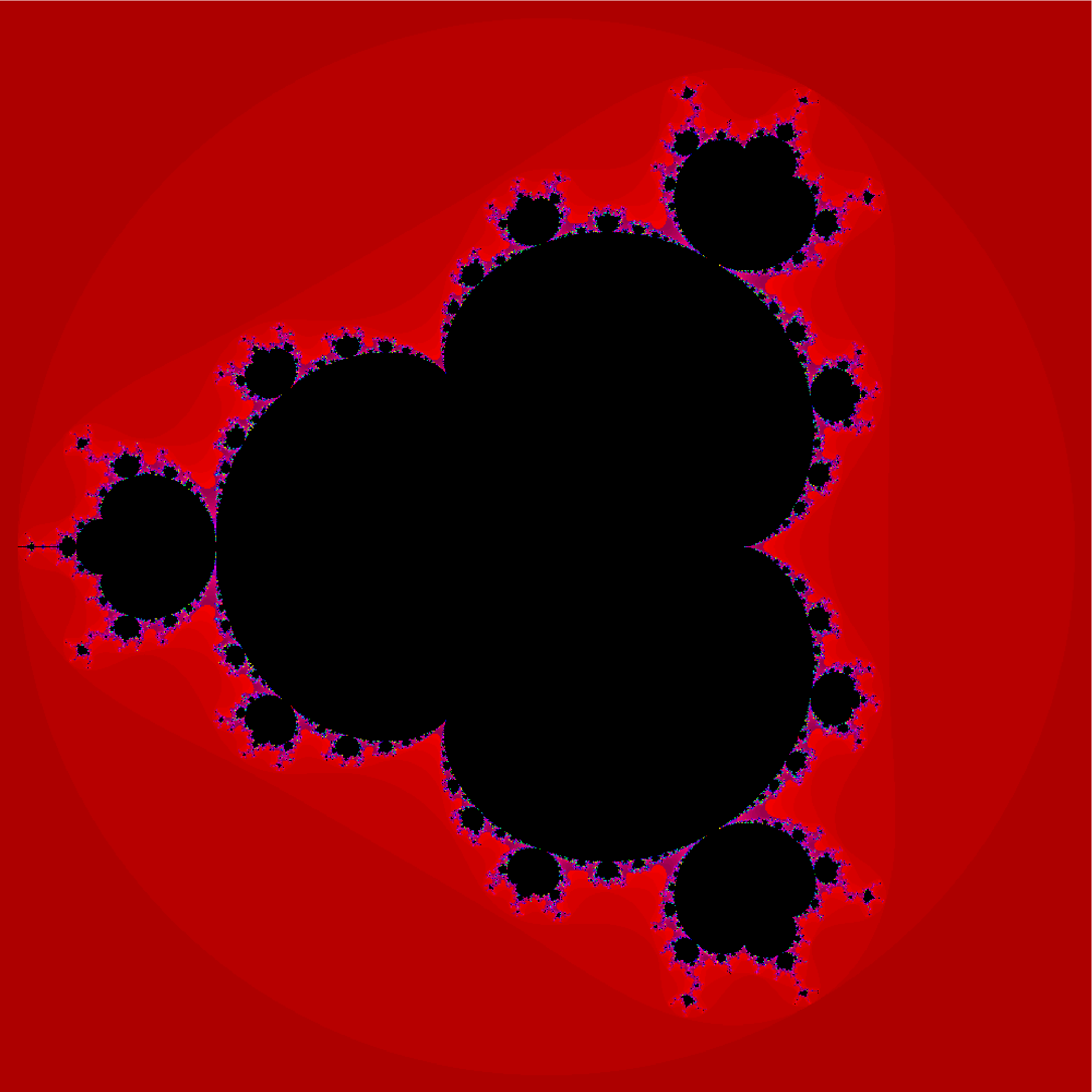}\label{fig2.2}}
\caption{Examples of $\mManp$ sets for $p=3,4$}
\end{figure}

We define the \textit{Multibrot} as follows:
\begin{definition}\label{d2.2.1}
Let $Q_{p,c}(z)=z^p+c$ a polynomial of degree $p\in \mN\setminus \left\lbrace  0,1\right\rbrace$. A \textit{Multibrot} set is the set of complex numbers $c$ which for all $m \in \mN$, the sequence $\Qpit$ is bounded, \textit{i.e.}
\begin{equation}
\mManp = \oa c \in \mC \, | \Qpit \text{ is bounded } \fa.
\end{equation}
\end{definition}
If we set $p=2$, we find the well-known Mandelbrot set. The next two theorems provide a method to visualize the $\mManp$ sets (see \cite{chine6} and \cite{Parise}).
\begin{theorem}\label{t2.2.1}
For all complex number $c$ in $\mManp$, we have $|c| \leq 2^{1/(p-1)}$.
\end{theorem}
To show Theorem \ref{t2.2.1}, we need the following lemma.
\begin{lemma}\label{l2.2.0.1}
Let $|c|^{p-1}>2$ with $p\geq 2$ an integer. Then, $|Q_{p,c}^m(0)|\geq |c|(|c|^{p-1}-1)^{m-1}$ for all natural number $m\geq 1$.
\end{lemma}
\begin{demo}
The proof is done by induction. For $m=1$, we have that $|Q_{p,c}(0)| = |c| = |c|(|c|^{p-1}-1)^{1-1}$. Suppose that the property is true for an integer $k\geq 1$. Then, for $k+1$, we obtain that
\begin{align*}
|Q_{p,c}^{k+1}(0)| &= |(Q_{p,c}^k(0))^p + c| \geq |Q_{p,c}^k(0)|^p-|c|\text{.}
\end{align*}
By the induction hypothesis and since $|c|^{p-1}>2$, we get the following inequalities
\begin{align*}
|Q_{p,c}^{k+1}(0)| &\geq |c|^p(|c|^{p-1}-1)^{p(k-1)} - |c|\\
&\geq |c|^p(|c|^{p-1}-1)^{k-1} - |c|(|c|^{p-1}-1)^{k-1}\\
&\geq |c|(|c|^{p-1}-1)^k\text{.}
\end{align*}
Hence, the property holds for $k+1$ and then it holds for all natural number $m\geq 1$.$\square$
\end{demo}
\begin{demo}[of Theorem \ref{t2.2.1}]
Suppose, by contradiction, there exists a complex number $c\in \mManp$ such that $|c|>2^{1/(p-1)}$. So, we have $|c|^{p-1}>2$. It follows from Lemma \ref{l2.2.0.1} that $|Q_{p,c}^m(0)|\geq |c|(|c|^{p-1}-1)^{m-1}$ for all $m\geq 1$. Then, as $m\rightarrow \infty$, $|Q_{p,c}^m(0)|$ tends to infinity since $|c|^{p-1}-1>1$. Thus, the sequence $\Qpit$ is unbounded, so $c\not \in \mManp$ by the definition of a Multibrot set.
\end{demo}
\begin{theorem}\label{t2.2.2}
A complex number $c$ is in $\mManp$ iff $|Q_{p,c}^m(0)|\leq 2^{1/(p-1)}$ $\forall m \in \mN$.
\end{theorem}
We need another lemma to prove Theorem \ref{t2.2.2}.
\begin{lemma}\label{l2.2.0.2}
Let $|c|\leq 2^{1/(p-1)}$ with $p\geq 2$ an integer. Suppose there exists an integer $n\geq 1$ such that $|Q_{p,c}^n(0)|$ $= 2^{1/(p-1)}+\delta$ with $\delta > 0$. Then, we have this following inequality: $|Q_{p,c}^{n+m}(0)| \geq 2^{1/(p-1)}+(2p)^m\delta$, $\forall m\geq 1$.
\end{lemma}
\begin{demo}
The proof is done by induction. For $m=1$, we have that $|Q_{p,c}^{n+1}(0)| \geq |Q_{p,c}^n(0)|^p - |c|$. By the hypothesis, we get $|Q_{p,c}^{n+1}(0)|\geq (2^{1/(p-1)}+\delta)^p - 2^{1/(p-1)}$. But, $(2^{1/(p-1)}+\delta)^p = \sum_{i=0}^p\binom{p}{i}2^{\frac{p-i}{p-1}}\delta^i$. Then,
\begin{align*}
(2^{1/(p-1)}+\delta)^p - 2^{1/(p-1)} &\geq 2^{p/(p-1)} + 2p\delta - 2^{1/(p-1)}\\
&=2^{1/(p-1)}(2^p-1) + 2p\delta\\
&\geq 2^{1/(p-1)}+2p\delta\text{.}
\end{align*}
Thus, $|Q_{p,c}^{n+1}(0)| \geq 2^{1/(p-1)}+2p\delta$ and the property holds for $m=1$. Now, suppose that the property is true for an integer $k\geq 1$. Then, for $k+1$, we have that $|Q_{p,c}^{n+k+1}(0)| \geq |Q_{p,c}^{n+k}(0)|^p - |c|$. Then, by the induction hypothesis, we get the following chain of inequalities
\begin{align*}
|Q_{p,c}^{n+k}(0)|^p - |c| &\geq (2^{1/(p-1)}+(2p)^k\delta)^p - |c|\\
&\geq 2^{p/(p-1)} + 2p(2p)^k\delta - 2^{1/(p-1)}\\
&\geq 2^{1/(p-1)} + (2p)^{k+1}\delta\text{.}
\end{align*}
Consequently, the property holds for $k+1$ and thus it holds for any natural number $m\geq 1$. 
\end{demo}
\begin{demo}\textbf{(of Theorem \ref{t2.2.2})}

\begin{enumerate}
\item[$\Rightarrow$)] Let $c \in \mManp$. By Theorem \ref{t2.2.1}, we know that $|c|\leq 2^{1/(p-1)}$. Suppose there exists a $n\geq 1$ such that $|Q_{p,c}^n(0)| > 2^{1/(p-1)}$, that is $|Q_{p,c}^n(0)| = 2^{1/(p-1)}+\delta$ with $\delta > 0$. Then, by Lemma \ref{l2.2.0.2}, we obtain that $|Q_{p,c}^{n+m}(0)| \geq 2^{1/(p-1)}+(2p)^m\delta$ for all $m\geq 1$. Then, $|Q_{p,c}^{n+m}(0)| \to \infty$ as $m$ tends to infinity since $2p>1$. Thus, $c\not\in \mManp$. This is a contradiction.
\item[$\Leftarrow$)] This is a direct consequence of the definition of Multibrot sets since the sequence $\Qpit$ is bounded by $2^{1/(p-1)}$ for all natural number $m\geq 1$.
\end{enumerate}
\end{demo}
Theorem \ref{t2.2.2} provides a criterion to decide whenever a complex number $c$ belongs to the set $\mManp$. The algorithm used to generate the figures is described in \cite{Gujar}. We use the preset limit $L=2^{1/(p-1)}$ and the magnitude of maximum iterations $M=1000$. The images are generated in a square of $1000 \times 1000$ pixels. Figures \ref{fig2.1}, \ref{fig2.2} illustrate, respectively, the sets $\mManp$ for the values $p=3$ and $p=4$.

Now, let $\mathfrak{M}$ denote the family of all generalized Mandelbrot sets $\mManp$, \textit{i.e.} $\mathfrak{M}:=\oa \mManp \, | \, p\geq 2\fa$. The family $\mathfrak{M}$ has the following property (see \cite{Chine1}).
\begin{theorem}\label{t2.2.3}
All member of the family $\mathfrak{M}$ is a connected set.
\end{theorem}

We have also some other properties related to the polynomial $Q_{p,c}(z)$ when we iterate it from $z_0=0$. The proofs can be found in \cite{Parise}.
\begin{lemma}\label{l2.2.1}
Set $c>0$ where $c$ is a real number. Then, the sequence $\Qpit$ is strictly ascendant. Furthermore, if the sequence $\Qpit$ is bounded, then it converges to $c_0>0$.
\end{lemma}
\begin{lemma}\label{l2.2.2}
Set $c<0$ where $c$ is a real number. Then, the sequence $\Qpit$ is strictly descendant when $p$ is odd. Furthermore, if the sequence $\Qpit$ is bounded, then it converges to $c_0<0$.
\end{lemma}

These properties will be useful to prove our next result for the intersection of $\mManth$ with the real axis.

\subsection{Roots of a third-degree polynomial}
Let $P(x)=x^3+bx^2+cx+d$ denote a monic cubic polynomial with real coefficients. Set $y=x+\frac{b}{3}$ as a Möbius transformation. It reduces the polynomial $P(x)$ to the polynomial $Q(y)=y^3+py+q$ where $p=c-\frac{b^2}{3}$ and $q=\frac{2b^3}{27}-\frac{cb}{3}+d$. In that case, searching for the roots of $P(x)$ is equivalent to search the roots of $Q(y)$.

A complex number $z$ is a root of $Q(y)$ iff there exist two complex numbers $y_1$ and $y_2$ such that $z=y_1+y_2$ and
\begin{align}
y_1^3+y_2^3+q=0\\ y_1y_2=-\dfrac{p}{3}
\label{syst1}
\end{align}
(see \cite{Parise}). The last equations can be rewritten as the following equivalent system
\begin{align}
y_1^3+y_2^3+q=0\notag \\ y_1^3y_2^3=-\dfrac{p^3}{27}
\label{syst2}\text{.}
\end{align}
With respect to \eqref{syst2}, we remark that $y_1^3$, $y_2^3$ are roots of the polynomial $T(t)=t^2-(y_1^3+y_2^3)t+y_1^3y_2^3=t^2+qt-\frac{p^3}{27}$ where its discriminant $\Delta$ is equal to $q^2+\frac{4p^3}{27}$. To settle information on the roots of $Q(y)$ and so, to get information from the roots of $P(x)$, we are interesting about the sign of $D=27\Delta =27q^2+4p^3$. In fact, one can prove the following result (see \cite{HBook}, \cite{EQU34} and \cite{Parise}).
\begin{theorem}\label{t2.1.1}
Let $P(x)=x^3+bx^2+cx+d$ where $b,c,d \in \mathbb{R}$. If $D=4c^3+27d^2+4db^3-b^2c^2-18bcd$ and 
\begin{itemize}
\item[i)] $D>0$, then $P(x)$ has one real root and two complex roots;
\item[ii)] $D=0$, then $P(x)$ has three real roots which one is of multiplicity 2;
\item[iii)] $D<0$, then $P(x)$ has three distinct real roots.
\end{itemize}
\end{theorem}
\begin{demo}
Developing the expression of $D = 27\Delta$ where $\Delta = q^2 + \frac{4p^3}{27}$ is the discriminant of $T(t)$ gives the expression of $D$ in the statement.

Now, we study the relation between the roots of $Q(y)$ and the sign of $D$.
\begin{enumerate}
\item[i)] If $D > 0$, then $T(t)$ has two real roots $t_1$ and $t_2$. According to \eqref{syst2}, $t_1 = y_1^3$ and $t_2 = y_2^3$. So, from the remark above, the roots of $Q(y)$ are
\begin{align*}
z_1& = y_1 + y_2 = \sqrt[3]{t_1} + \sqrt[3]{t_2}\\
z_2& = \omega y_1 + \overline{\omega}y_2 = \omega \sqrt[3]{t_1} + \overline{\omega}\sqrt[3]{t_2}\\
z_3& = \overline{\omega}y_1 + \omega y_2 = \overline{\omega}\sqrt[3]{t_1} + \omega \sqrt[3]{t_2}
\end{align*}
where $\omega:=\frac{-1+\bi \sqrt{3}}{2}$ is a cubic root of the unity. This gives one real root and two complex roots.
\item[ii)] If $D = 0$, then $T(t)$ has a root $t = -\frac{q}{2}$ of multiplicity $2$. According to \eqref{syst2}, $y^3 = t$ where $y^3:=y_1^3 = y_2^3$. Then, the roots of $Q(y)$ are
\begin{align*}
z_1 &= 2y = 2\sqrt[3]{t}\\
z_2 &= z_3 = \omega y + \overline{\omega} y = -\sqrt[3]{t}\text{.}
\end{align*}
These last roots are all real, and one of them (the root $z_2$) is of multiplicity 2.
\item[iii)] If $D < 0$, then the roots $t$ and $\overline{t}$ are complex roots of $T(t)$. Let $t = re^{\bi\theta }$ and set $y:= t^{1/3}$ a cubic root of $t$. Then, the roots of $Q(y)$ are
\begin{align*}
z_1& = y + \overline{y} = t^{1/3} + \overline{t}^{1/3}\\
z_2& = \omega y + \overline{\omega}\, \overline{y} = \omega t^{1/3} + \overline{\omega}\,\overline{t}^{1/3}\\
z_3&= \overline{\omega} y + \omega\overline{y} = \overline{\omega}t^{1/3} + \omega \overline{t}^{1/3}\text{.}
\end{align*}
Since the three roots are the sum of a complex number with its conjugate, $z_1$, $z_2$ and $z_3$ are distinct real roots. This complete the proof.
\end{enumerate}

\end{demo}
\subsection{$\mManth$ set: the \textit{Mandelbric}}
In this subsection, we prove that the \textit{Mandelbric} set crosses the real axis on the interval $\left[-\frac{2}{3\sqrt{3}}, \frac{2}{3\sqrt{3}} \right]$ (see Theorem \ref{t2.3.1}). Before we go through the proof of Theorem \ref{t2.3.1}, we first establish some symmetries (see \cite{Lau} and \cite{Sheng}) in the \textit{Mandelbric}.
\begin{lemma}\label{l2.3.1}
Let $c \in \mManth$. Then $\cc \in \mManth$.
\end{lemma}
\begin{demo}
 Suppose $c \in \mManth$. Then, by Theorem \ref{t2.2.2}, $|Q_{3,c}^m(0)| \leq \sqrt{2}$ $\forall m \in \mN$. Since $Q_{3,\cc}^m(0)=\overline{Q_{3,c}^m(0)}$, $|Q_{3,\cc}^m(0)|=|Q_{3,c}^m(0)|$ for all $m \in \mN$. Thus, $\cc \in \mManth$.
\end{demo}
Lemma \ref{l2.3.1} provides that $\mManth$ is symmetrical about the real axis. The next lemma is a step forward to show that $\mManth$ is symmetrical about the imaginary axis.
\begin{lemma}\label{l2.3.2}
Let $c=x+ y\bi$ where $x,y \in \mR$. If $c \in \mManth$, then $-c \in \mManth$.
\end{lemma}
\begin{demo}
Let $c=x+  y\bi$ where $x,y \in \mR$. If $c \in \mManth$, then by Theorem \ref{t2.2.2}, $|Q_{3,c}^m(0)|\leq \sqrt{2}$ $\forall m \in \mN$. By induction, we remark that $\forall m \in \mN$, $Q_{3,-c}^m(0)=-Q_{3,c}^m(0)$, and so $|Q_{3,-c}^m(0)|=|Q_{3,c}^m(0)|\leq \sqrt{2}$ $\forall m \in \mN$. Thus, $-c \in \mManth$.
\end{demo}
\begin{corollary}\label{c2.3.1}
Let $c=x+ y\bi$ where $x,y \in \mR$. If $c \in \mManth$, then $c'=-x+y\bi $ is in $\mManth$.
\end{corollary}
\begin{demo}
Let $c=x+ y\bi$ and $c \in \mManth$. We want to prove that $c'=-x+ y\bi\in\mManth$. By hypothesis and Lemma \ref{l2.3.1}, $\cc \in \mManth$. Therefore, by Lemma \ref{l2.3.2}, $-\cc \in \mManth$. Since $-\cc=-x+ y\bi$, we have that $c'\in \mManth$.
\end{demo}
With this last result, the next proof will be simplifies.
\begin{theorem}\label{t2.3.1}
The Mandelbric crosses the real axis on the interval $\left[ \frac{-2}{3\sqrt{3}},\, \frac{2}{3\sqrt{3}} \right]$.
\end{theorem}
\begin{demo}
By the Corollary \ref{c2.3.1}, we can restrict our proof to the interval $\left[ 0,\, \frac{2}{3\sqrt{3}} \right]$. Let $R_{3,c}(x)=x^3-x+c$ where $c \in \mR$ and $D=-4+27c^2$. We start by showing that no point $c>\frac{2}{3\sqrt{3}}$ lies in $\mManth$. In this case, $D>0$ and $R_{3,c}$ has a real root (see Theorem \ref{t2.1.1}), and this root is given by
\begin{equation}
x_0=\sqrt[3]{ -\dfrac{c}{2} + \dfrac{\sqrt{c^2-\frac{4}{27}}}{2}} + \sqrt[3]{ -\dfrac{c}{2} - \dfrac{\sqrt{c^2-\frac{4}{27}}}{2} }.
\end{equation}
Suppose that $c \in \mManth$, \textit{i.e.} $\Qit$ is bounded. Then, Lemma \ref{l2.2.1} implies that $\Qit$ is strictly ascendant and it converges to $c_0>0$. Since $Q_{3,c}^m(0)$ is a polynomial function for all $m\in \mN$, we have that
\begin{align}
c_0&=\lim_{m \rightarrow \infty} Q_{3,c}^{m+1}(0)\\
&= Q_{3,c}( \lim_{m \rightarrow \infty} Q_{3,c}^{m}(0) )\\
&= Q_{3,c}(c_0)\text{.}
\end{align}
Thus, $c_0$ is a real root of $R_{3,c}$ and $c_0=x_0$. However, since $\frac{c}{2}>-\frac{c}{2}$, we have that
\begin{equation}
x_0= \sqrt[3]{ -\dfrac{c}{2} + \dfrac{\sqrt{c^2-\frac{4}{27}}}{2} } + \sqrt[3]{ -\dfrac{c}{2} - \dfrac{\sqrt{c^2-\frac{4}{27}}}{2} }<0.
\end{equation}
This is a contradiction with $x_0=c_0>0$. Thus, $c \not \in \mManth$.

Next, we show that for $0 \leq c \leq \frac{2}{3\sqrt{3}}$, $c$ lies in $\mManth$. Obviously, $c=0$ is in $\mManth$. Suppose that $0 < c \leq \frac{2}{3\sqrt{3}}$. In this case, $D \leq 0$ and by the Theorem \ref{t2.1.1}, $R_{3,c}(x)$ has the following three real roots (see \cite{Parise}):
\begin{equation}\label{e2.3.1}
{\left( -\dfrac{c}{2} + \bi\dfrac{\sqrt{-c^2+\frac{4}{27}}}{2} \right)}^{1/3} + {\left( -\dfrac{c}{2} -\bi \dfrac{\sqrt{-c^2+\frac{4}{27}}}{2} \right)}^{1/3}.
\end{equation}
Following De Moivre's formula, one of these roots can be expressed as follows:
\begin{equation}
a=\dfrac{2}{\sqrt{3}}\cos \left( \frac{\theta}{3}\right)
\end{equation}
for $c\in \left( 0, \, \frac{2}{3\sqrt{3}} \right]$ and $\theta=\arctan \left(\frac{\sqrt{-D}}{-3c\sqrt{3}}\right)+\pi$ where $\pi \leq  \theta < \frac{\pi}{2}$. We prove by induction that $|Q_{3,c}^m(0)|<a$ $\forall m \in \mN$. For $m=1$, we have that $|Q_{3,c}(0)|=|c|<a$ because $|c| < \frac{1}{\sqrt{3}} \leq a$. Indeed, since $\pi \leq  \theta < \frac{\pi}{2}$, we obtain $\frac{1}{\sqrt{3}} \leq a < 1$. Now, suppose that $|Q_{3,c}^k(0)| < a$ for a $k \in \mN$. Then, since $R_{3,c}(a)=a^3-a+c=0$ and $c>0$,
\begin{equation*}
|Q_{3,c}^{k+1}(0)|=|(Q_{3,c}^k(0))^3+c|\leq |(Q_{3,c}^k(0))|^3+|c| < a^3+c=a.
\end{equation*}
Thus, the proposition is true for $k+1$ and $|Q_{3,c}^m(0)|<a$ $\forall m \in \mN$. Since $a \leq \sqrt{2}$, then by the Theorem \ref{t2.2.2} we have $c \in \mManth$ .

In conclusion, $\mManth \cap \mR_+ = \left[ 0, \, \frac{2}{3\sqrt{3}} \right]$.
\end{demo}

\section{\textit{Multibrot} sets for hyperbolic numbers}
Previously, we treated the \textit{Multibrot} sets for complex numbers. In this section, we propose an extension of the Mandelbrot set for hyperbolic numbers called the \textit{Hyperbrots} and we prove that $\mManth$ for hyperbolic numbers denoted by $\mHth$ is a square of side length $\frac{2}{3\sqrt{3}}\sqrt{2}$.

\subsection{Definition of the sets $\mathcal{H}^p$}
Based on the works of Metzler and Senn (see \cite{MET} and \cite{Senn} respectively) on the hyperbolic Mandelbrot set, we define the Hyperbrots as follows:
\begin{definition}\label{d3.1}
Let $Q_{p,c}(z)=z^p+c$ where $z,c \in \mD$ and $p\geq 2$ an integer. The Hyperbrots are defined as the sets
\begin{equation}
\mHyb^p:= \oa c \in \mD \, | \, \Qpit \text{ is bounded }  \fa \text{.}
\end{equation}
\end{definition}
Metzler proved that $\mHt$ is a square with diagonal length 2$\frac{1}{4}$ and of side length $\frac{9}{8} \sqrt{2}$. We use the same approach to prove that $\mHth$ is a square with diagonal length $\frac{4}{3\sqrt{3}}$ and with side length $\frac{2}{3\sqrt{3}}\sqrt{2}$.
For the next part, we note a hyperbolic numbers $z$ as $(u,v)^{\top}$ and the fixed number $c$ as $(a,b)^{\top}$ where $\top$ is the transpose of a column vector in $\mR^2$.

\subsection{Special case: $\mathcal{H}^3$}
First, we recall some of the basic tools introduced by Metzler.
\begin{definition}\label{d3.2}
Let $(u,v)^{\top}, (x,y)^{\top} \in \mR^2$. We define two multiplication operations $\diamond$ and $\ast$ on $\mR^2$ as
\begin{equation}\label{e3.1}
\begin{pmatrix}
u\\v
\end{pmatrix} \diamond
\begin{pmatrix}
x\\y
\end{pmatrix}:=
\begin{pmatrix}
ux+vy\\vx+uy
\end{pmatrix}
\end{equation}
\begin{equation}\label{e3.2}
\begin{pmatrix}
u\\v
\end{pmatrix} \ast
\begin{pmatrix}
x\\y
\end{pmatrix}:=
\begin{pmatrix}
ux\\vy
\end{pmatrix}.
\end{equation}
\end{definition}
\begin{remark}\label{r3.1}
The operation $\diamond$ corresponds to the multiplication $\cdot$ of two hyperbolic numbers as we adopted the two-dimensional vector notation. We use the symbols $\diamond_{\circ n}$ and $\ast_{\circ n}$ to denote the $n$ consecutive multiplications $\diamond \circ \diamond \circ \ldots \circ \diamond$ and $\ast \circ \ast \circ \ast \ldots \circ \ast$ respectively.
\end{remark}
With the usual addition operation $+$ on $\mR^2$, $(\mR^2,+,\diamond )$ and $(\mR^2,+,\ast )$ are commutative rings with unity.

\begin{figure}
\centering
\includegraphics[scale=0.25]{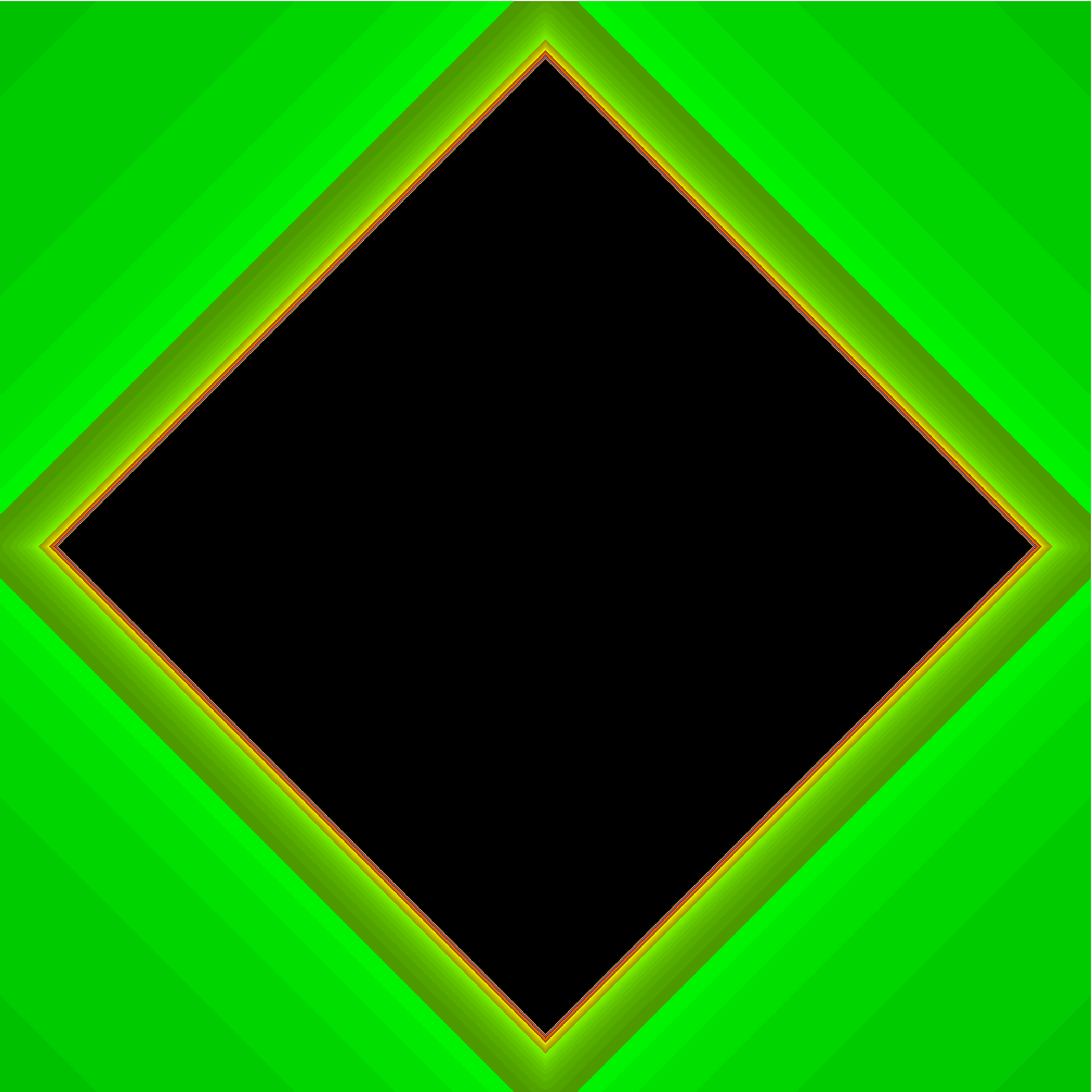}
\caption{\textit{Hyperbrot} $\mHth$: $-0.40\leq \pre (z) \leq 0.40$, $-0.40 \leq \pim (z) \leq 0.40$}\label{fig3.1}
\end{figure}

Let $T:\, \mR^2 \rightarrow \mR^2$ be the following matrix
\begin{equation}\label{e3.3}
T:=
\begin{pmatrix}
1 & -1 \\ 1 & 1
\end{pmatrix}\text{.}
\end{equation}
Then, $T$ is an isomorphism between $(\mR^2,+,\diamond )$ and $(\mR^2,+,\ast )$. Now, we define
\begin{equation}\label{e3.4}
\bHpc \vecxy := \vecxy \diamond_{\circ p} \vecxy + \vecab
\end{equation}
and we generalize a result that is included in the proof of Metzler for the case $p=2$.
\begin{lemma}\label{l3.1}
For all $m \in \mN$, we have that
\begin{equation}\label{e3.5}
T\bHpc^m \vecxy = \begin{pmatrix}
Q_{p,a-b}^m(x-y) \\ Q_{p,a+b}^m(x+y)
\end{pmatrix}
\end{equation}
where $T$ is the matrix of equation \eqref{e3.3} and $Q_{p,c}(z)=z^p+c$ with $z,c\in \mR$.
\end{lemma}
The proof can be found in \cite{Parise}. It is similar to the one of Metzler gave in his article (see \cite{MET}). We just replace $\diamond$ by $\diamond_{\circ p}$, $\ast$ by $\ast_{\circ p}$ and the degree $2$ of $P_{(a,b)}$ by the integer $p \geq 2$. Hence, Lemma \ref{l3.1} allows to separate the dynamics of $\bHpc \vecxy$ in terms of the dynamics of two real polynomials $Q_{p,a-b}(x-y)$ and $Q_{p,a+b}(x+y)$. We now use Lemma \ref{l3.1} and Theorem \ref{t2.3.1} to prove the next result.
\begin{theorem}\label{t3.2}
$\mHth = \oa c=(a,b)^{\top} \in \mR^2 \, | \, |a|+|b| \leq \frac{2}{3\sqrt{3}} \fa$.
\end{theorem}
\begin{demo}
By Lemma \ref{l3.1} and the remark right after, $\oa \bHpc^m (\mathbf{0}) \fa _{m=1}^{\infty}$ is bounded  iff the real sequences $\oa Q_{3,a-b}^m(0) \fa _{m=1}^{\infty}$ and $\oa Q_{3,a+b}^m(0) \fa _{m=1}^{\infty}$ are bounded . But, according to Theorem \ref{t2.3.1}, these sequences are bounded iff
\begin{align}
|a-b| \leq \frac{2}{3\sqrt{3}} \mbox{ and } |a+b|\leq \frac{2}{3\sqrt{3}} \label{e3.7}\text{.}
\end{align}
Then, by a simple computation we obtain $|a|+|b| \leq \frac{2}{3\sqrt{3}}$. Conversely, if $|a|+|b| \leq \frac{2}{3\sqrt{3}}$ is true, then by the properties of the absolute value, we obtain the inequalities in \eqref{e3.7}. Thus, we obtain the following characterization for $\mHth$, $\mHth = \oa c=(a,b)^{\top} \in \mR^2 \, | \, |a|+|b| \leq \frac{2}{3\sqrt{3}} \fa \text{.}$
\end{demo}
Figure \ref{fig3.1} represent an image of the set $\mHth$.

\section{Generalized Mandelbrot sets for tricomplex numbers}\label{sec5}
In this section, we use the set of tricomplex numbers to generalize the \textit{Multibrot} sets. Particularly, we give some basic properties of the tricomplex \textit{Multibrot} and we continue the exploration of the \textit{Multibrot} started in \cite{GarantRochon}. We will concentrate our exploration on the case $\mMan_3^3$ corresponding to the polynomial $Q_{3,c}(\eta)=\eta^3+c$ with $\eta,c \in \mM (3)$.

\subsection{Definition and properties of $\mathcal{M}_3^p$}
The authors of \cite{Chine1} defined the bicomplex \textit{Multibrot} sets as follows:
\begin{definition}\label{def5.1}
Let $Q_{p,c}=\zeta^p+c$ where $\zeta , c \in \mM (2)$ and $p \geq 2$ is an integer. The bicomplex \textit{Multibrot} set is defined as the set
\begin{equation*}
\mMan_2^p:= \oa c\in \mM (2) \, | \, \Qpit \text{ is bounded } \fa \text{.}
\end{equation*}
\end{definition}
In \cite{Chine1}, they proved that $\mMan_2^p$ can be expressed as a $\mM (2)$-\textit{cartesian} set and is connected. In the next theorem, we improve their result concerning the bounded discus of $\mMan_2^p$ in conformity with our Theorem \ref{t2.2.1}.
\begin{theorem}\label{t5.1}
Let $\mMan_2^p$ denote the bicomplex \textit{Multibrot} sets for integers $p\geq 2$. Then, the following inclusions hold:
\begin{equation}
\mMan_2^p \subset \Ol{\mathbf{D_2}(0,2^{\frac{1}{p-1}},2^{\frac{1}{p-1}})} \subset \Ol{\mathbf{B_2^1}(0,2^{\frac{1}{p-1}})}\text{.}
\end{equation}
\end{theorem}
\begin{demo}
We know from \cite{Chine1} that $\mMan_2^p=\mMan_1^p \times_{\eb}\mMan_1^p$. Moreover, by Theorem \ref{t2.2.1}, we know that $\mMan_1^p \subset \Ol{\mathbf{B_1^1}(0,2^{\frac{1}{p-1}})}$. So, combining the both last statements, we proved the left inclusion. For the right inclusion,
we use this following inclusion $\Ol{\mathbf{D_2}(a,\bf{r_1},\bf{r_2})} \subset \Ol{\mathbf{B_2^1}(a,\sqrt{\frac{\bf{r_1}^2+\bf{r_2}^2}{2}})}$ (see \cite{Baley}) with $a=0$ and $\bf{r_1}=\bf{r_2} = 2^{\frac{1}{p-1}}$.
\end{demo}

Now, the tricomplex \textit{Multibrot} sets are defined analogously to the bicomplex ones:
\begin{definition}\label{d5.2}
Let $Q_{p,c}=\eta^p+c$ where $\eta , c \in \mM (3)$ and $p\geq 2$ an integer. The tricomplex \textit{Multibrot} set is defined as the set
\begin{equation}
\mMan_3^p:=\oa c \in \mM (3) \, | \, \Qpit \text{ is bounded } \fa \text{.}
\end{equation}
\end{definition}
We have the following theorem that characterizes $\mMan_3^p$ set as a $\mM (3)$-\textit{cartesian} product of $\mMan_2^p$.
\begin{theorem}\label{t5.2}
$\mMan_3^p=\mMan_2^p \times_{\ett} \mMan_2^p$.
\end{theorem}
\begin{demo}
Let $c=c_1 + c_2\bb=(c_1- c_2\bt)\ett + (c_1+ c_2\bt)\etc$ as a tricomplex numbers. So, by Definition \ref{d5.2}, $c \in \mMan_3^p$ iff $\oa Q_{p,c}^m(0) \fa$ is bounded. However, from Theorem \ref{theo2.2}, $Q_{p,c}^m(0)$ can be expressed with the idempotent representation as follows:
\begin{equation}
Q_{p,c}^m(0)=Q_{p,c_1- c_2\bt}^m(0)\ett + Q_{p,c_1+ c_2\bt}^m(0)\etc
\end{equation}
for all $m \in \mN$. Moreover, in \cite{Baley}, it is proved for the general case of multicomplex numbers that:
\begin{equation}
\Vert \zeta \Vert_n=\sqrt{\dfrac{\Vert \zeta_1 - \zeta_2 \bnm \Vert_{n-1}^2+ \Vert \zeta_1 + \zeta_2 \bnm \Vert_{n-1}^2}{2}}\label{e5.1}
\end{equation}
where $\zeta=\zeta_1+ \zeta_2 \bn\in \mM (n)$.
So, $\oa Q_{p,c}^m(0) \fa_{m=1}^{\infty}$ is bounded iff $\oa Q_{p,c_1- c_2\bt}^m(0) \fa_{m=1}^{\infty}$ and $\oa Q_{p,c_1+ c_2\bt}^m(0)\fa_{m=1}^{\infty}$ are bounded. By Definition \ref{def5.1}, we obtain that $c_1- c_2\bt, c_1+ c_2\bt \in \mMan_2^p$. Thus, $c=(c_1- c_2\bt)\ett + (c_1+ c_2\bt)\etc \in \mMan_2^p \times_{\ett} \mMan^p_2$.
\end{demo}
If we combine Theorem \ref{t5.2} with Theorem \ref{t2.2.1}, we get the following statement.
\begin{theorem}\label{t5.3}
Let $\mMan_3^p$ be the tricomplex \textit{Multibrot} set for $p\in\mN\backslash\{0,1\}$. Then the following inclusion holds:
\begin{equation}
\mMan_3^p \subset \Ol{\bf{D_3}}(0,2^{\frac{1}{p-1}},2^{\frac{1}{p-1}})\text{.}
\end{equation}
\end{theorem}

Finally, in \cite{Chine1}, it is proved that the sets $\mMan_2^p$ is connected $\forall p\in\mN\backslash\{0,1\}$. We obtain the same property for $\mMan_3^p$.
\begin{theorem}\label{t5.4}
$\mMan_3^p$ is a connected set.
\end{theorem}
\begin{demo}
Let define the function $\Gamma_2:\, X_1 \times X_2 \rightarrow X_1 \times_{\ett}X_2$ with $X_1,X_2 \subset \mM (2)$ and $X=X_1 \times X_2\subset \mM (3)$ by $\Gamma_2(u_1,u_2)=u_1\ett + u_2\etc$. Obviously, $\Gamma_2$ is an homeomorphism. So, if $X_1,X_2$ are connected sets, then $X$ is also a connected set. Thus, by Theorem \ref{d5.2}, $\mMan_3^p=\mMan_2^p \times_{\ett}\mMan_3^p$ and since $\mMan_2^p$ is connected (see \cite{Chine1}), $\mMan_3^p$ is also a connected set $\forall p\in\mN\backslash\{0,1\}$.
\end{demo}

Theorem \ref{t5.3} is useful to generate the divergence layers of the tricomplex \textit{Multibrot} sets. We use this information to draw the images of the next part. Moreover, we conjecture that the Fatou-Julia Theorem is true for tricomplex \textit{Multibrot} sets and use it to give more information about the topology of the sets. For a statement of the generalized Fatou-Julia Theorem, we refer the reader to \cite{GarantRochon}.

\subsection{Principal 3D slices of the set $\mathcal{M}_3^3$}
We want now to visualize the tricomplex \textit{Multibrot} sets. Since there are in eight dimensions, we take the same approach from \cite{GarantRochon} to accomplish this goal. In that way, we may denote the principal 3D slice for a specific tricomplex \textit{Multibrot} set as $\mTet^p$ and define it as the set
\begin{equation}\label{eq5.2.1}
\mTet^p:= \mTet^p(\bim , \bk , \bil ) =\oa c \in \mT (\bim , \bk , \bil ) \, | \, \oa Q_{p,c}^m(0) \fa_{m=1}^{\infty} \text{ is bounded } \fa
\text{.}
\end{equation}
So the number $c$ has three of its components that are not equal to zero. In total, there are 56 possible combinations of principal 3D slices. To attempt a classification of these slices, we introduce a relation $\sim$ (see \cite{GarantRochon}).
\begin{definition}\label{def5.2.1}
Let $\mTet_1^p( \bim , \bk , \bil )$ and $\mTet_2^p( \bn , \biq , \bis )$ be two 3D slices of a tricomplex \textit{Multibrot} set $\mathcal{M}_3^p$ that correspond, respectively, to $Q_{p,c_1}$ and $Q_{p,c_2}$. Then, $\mTet_1^p \sim \mTet_2^p$ if there exists a bijective linear function $\varphi : \mathrm{span}_{\mR}\oa 1,\bim , \bk , \bil \fa \rightarrow \mathrm{span}_{\mR} \oa 1, \bn , \biq , \bis \fa$ such that $(\varphi \circ Q_{p,c_1} \circ \varphi^{-1})(\eta )=Q_{p,c_2}(\eta )$ $\forall \eta \in \mathrm{span}_{\mR}\oa 1,\bn , \biq , \bis \fa$. In that case, we say that $\mTet_1^p$ and $\mTet_2^p$ have the same dynamics.
\end{definition}
If two 3D slices are in relationship in term of $\sim$, then we also say that they are symmetrical. This comes from the fact that their visualizations by a computer give the same images. In \cite{Parise}, it is showed that $\sim$ is also a equivalent relation on the set of 3D slices of $\mMan_3^p$. For the rest of this article, we focus on the principal slices of the $\mMan_3^3$ set, also called the tricomplex \textit{Mandelbric} set.

\begin{figure}
\centering
\subfigure[$\mathcal{T}^3(1,\mathbf{i_1},\mathbf{i_2})$: \textit{Tetrabric}]{%
\includegraphics[scale=0.3]{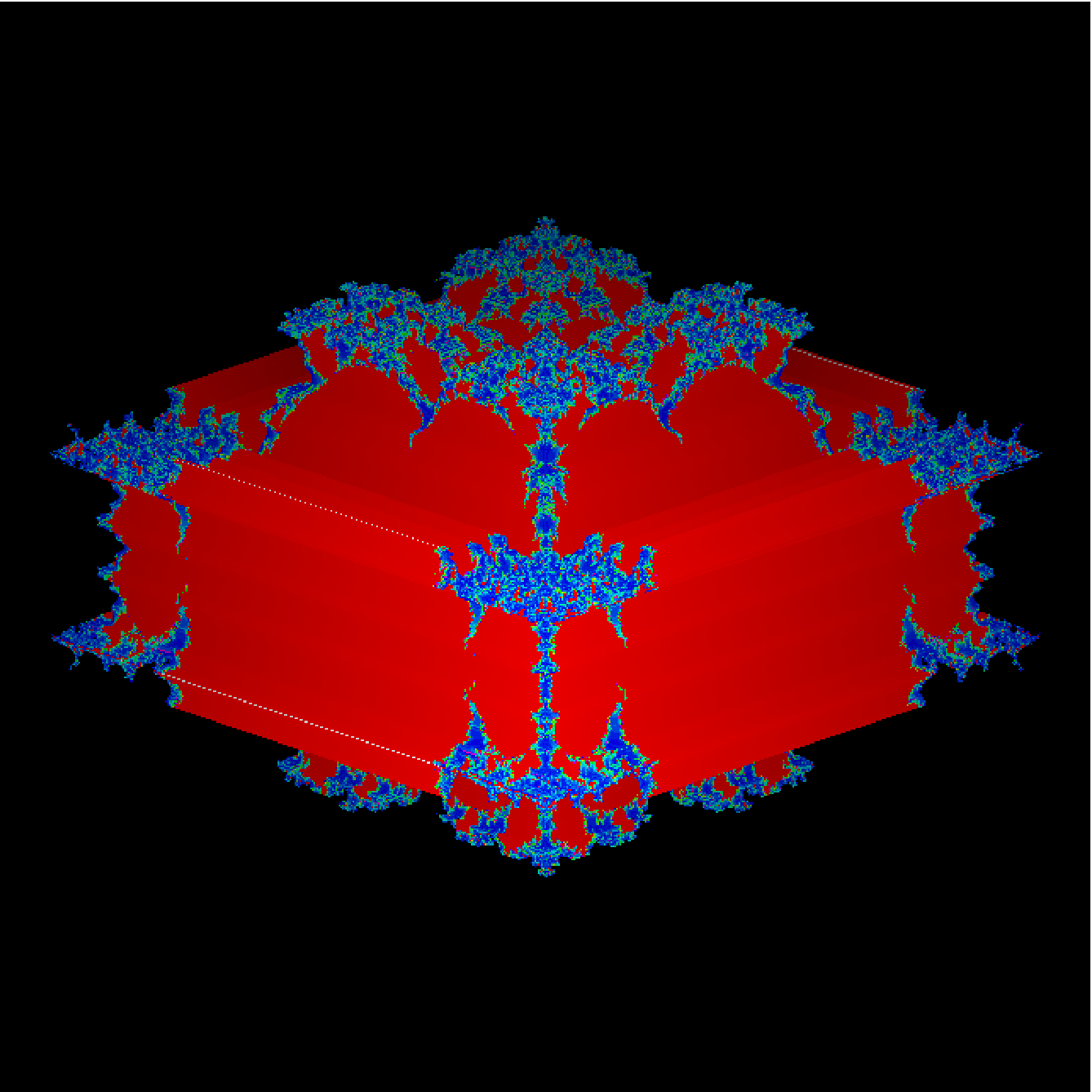}\label{fig5.1}}
\subfigure[$\mathcal{T}^3(1,\mathbf{j_1},\mathbf{j_2})$: \textit{Perplexbric}]{%
\includegraphics[scale=0.3]{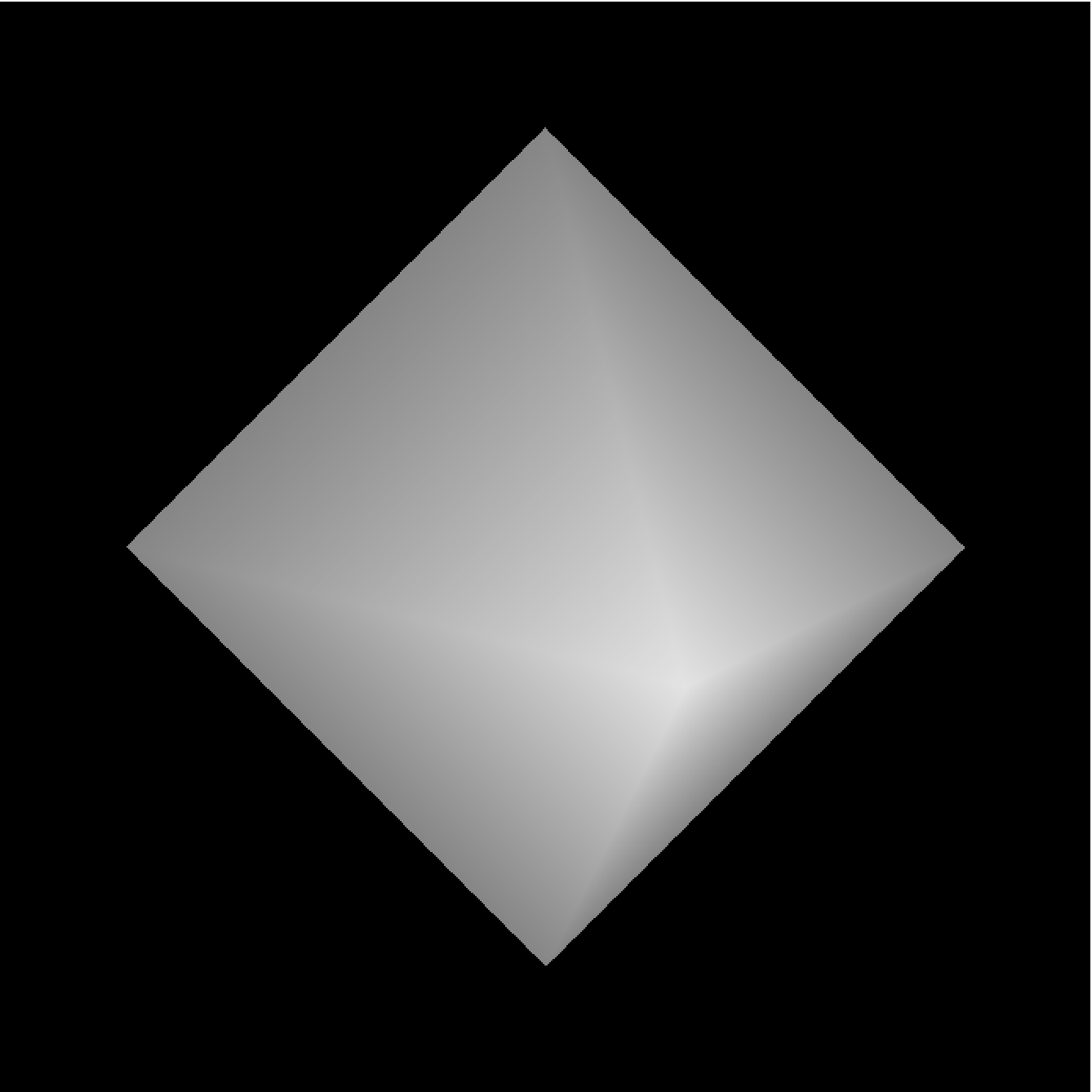}\label{fig5.2}}
\subfigure[$\mathcal{T}^3(\mathbf{i_1},\mathbf{i_2},\mathbf{j_1})$]{%
\includegraphics[scale=0.3]{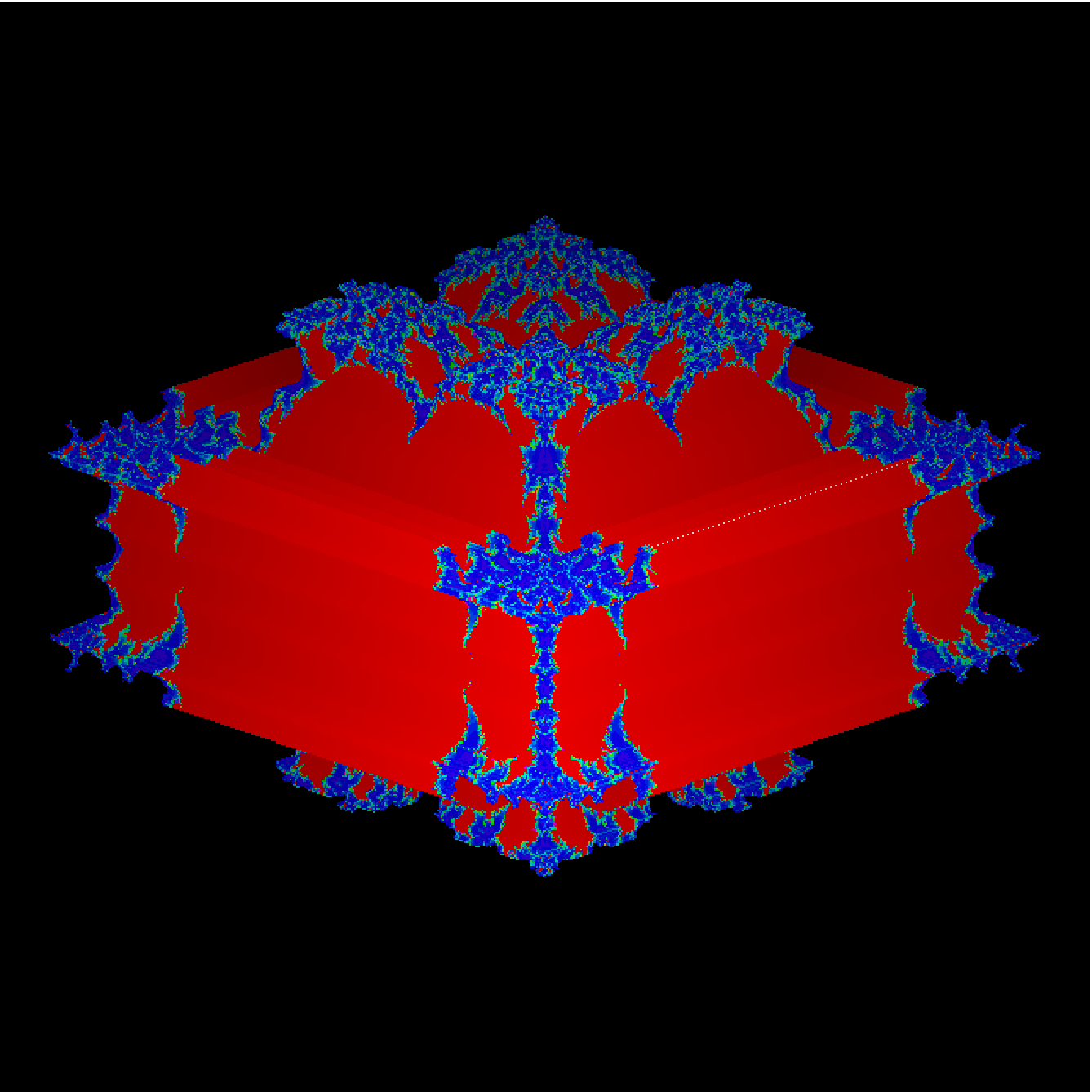}\label{fig5.3}}
\subfigure[$\mathcal{T}^3(1,\mathbf{i_1},\mathbf{j_1})$]{%
\includegraphics[scale=0.3]{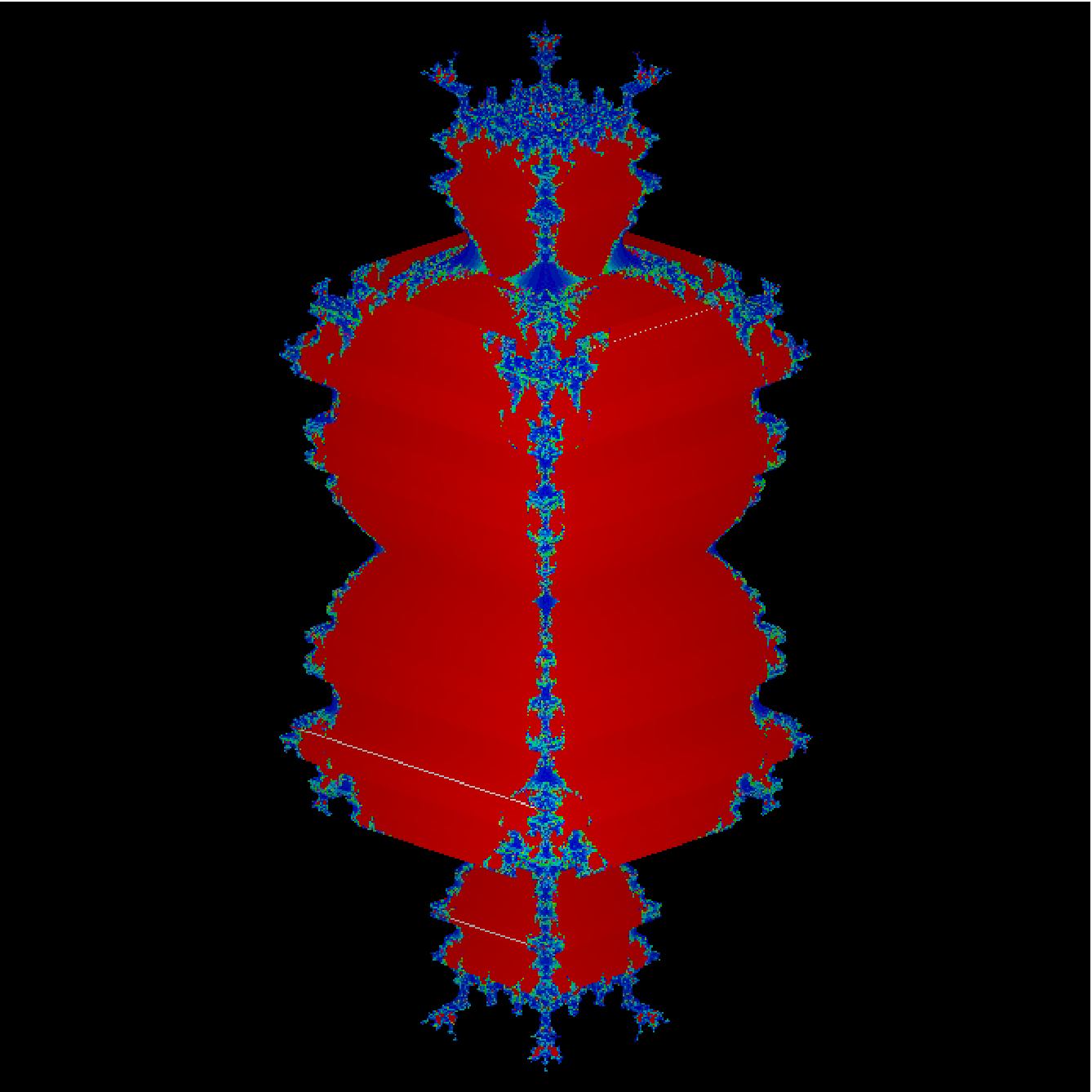}\label{fig5.4}}
\caption{Four 3D slices of the \textit{Mandelbric}}\label{figAllG1}
\end{figure}

\begin{figure}
\centering
\subfigure[$\mathcal{T}^3(\mathbf{i_1},\mathbf{i_2},\mathbf{i_3})$]{%
\includegraphics[scale=0.3]{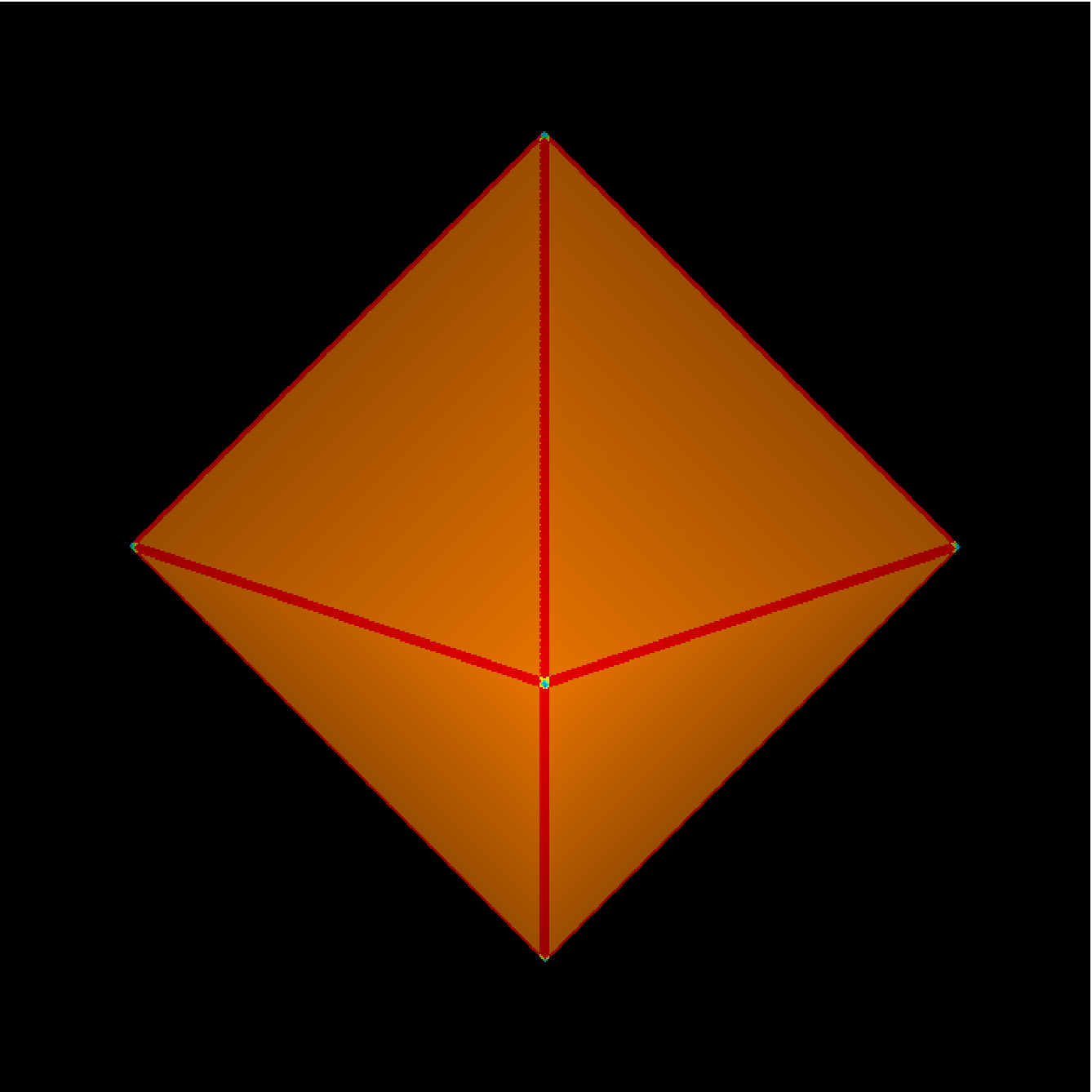}\label{fig5.5}}
\subfigure[$\mathcal{T}^3(\mathbf{i_1},\mathbf{i_2},\mathbf{j_2})$:]{%
\includegraphics[scale=0.3]{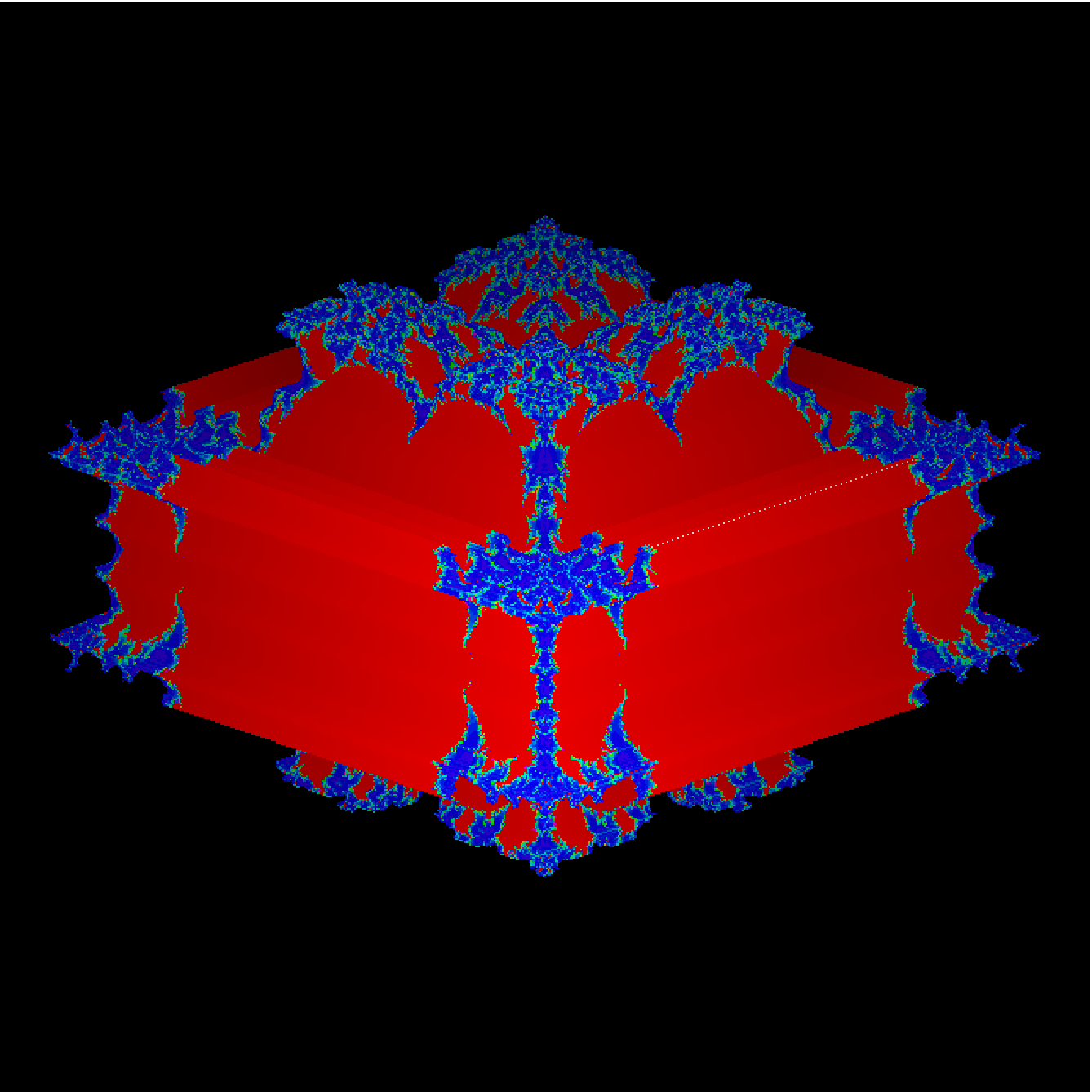}\label{fig5.6}}
\subfigure[$\mathcal{T}^3(\mathbf{i_1},\mathbf{j_1},\mathbf{j_2})$]{%
\includegraphics[scale=0.3]{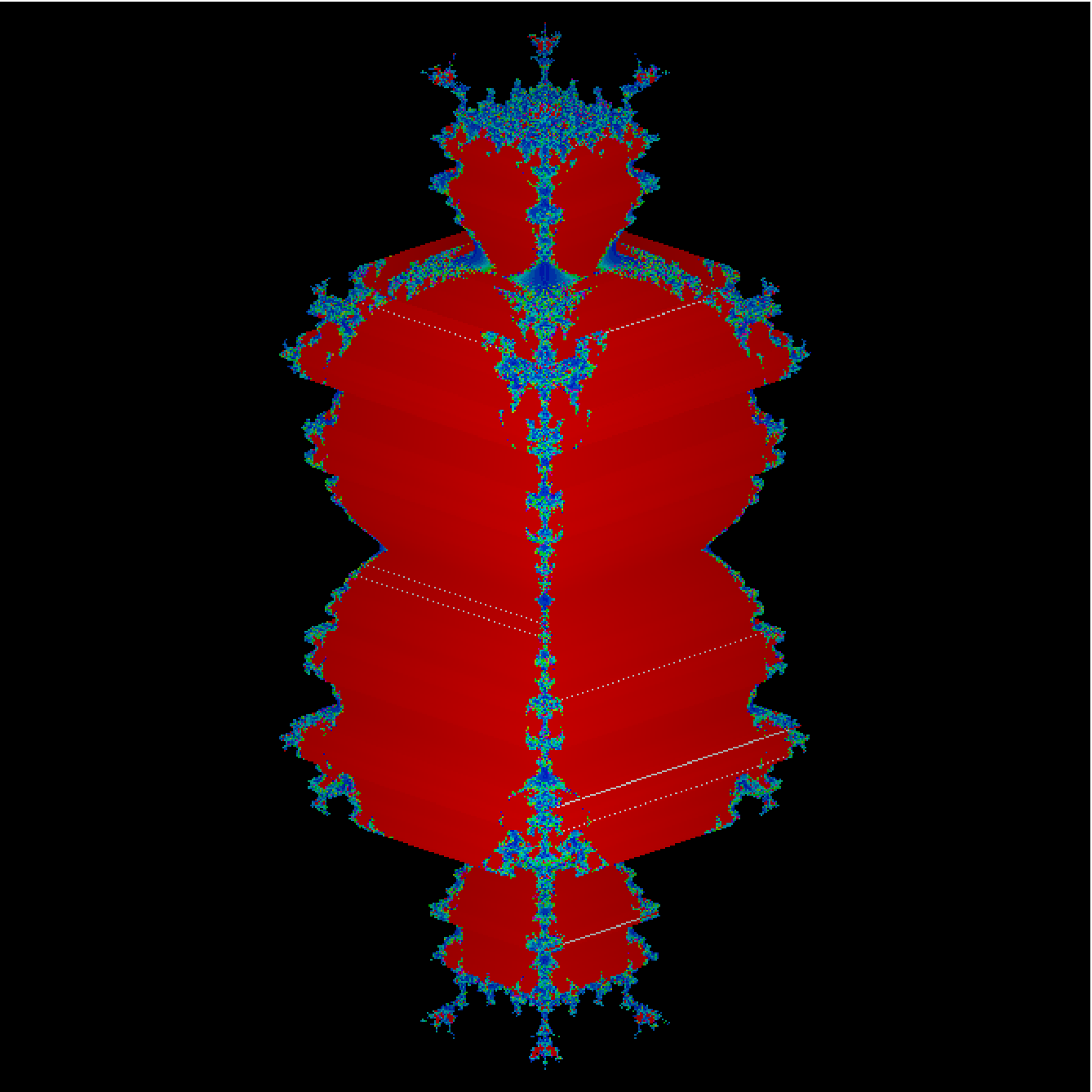}\label{fig5.7}}
\subfigure[$\mathcal{T}^3(\mathbf{j_1},\mathbf{j_2},\mathbf{j_3})$]{%
\includegraphics[scale=0.3]{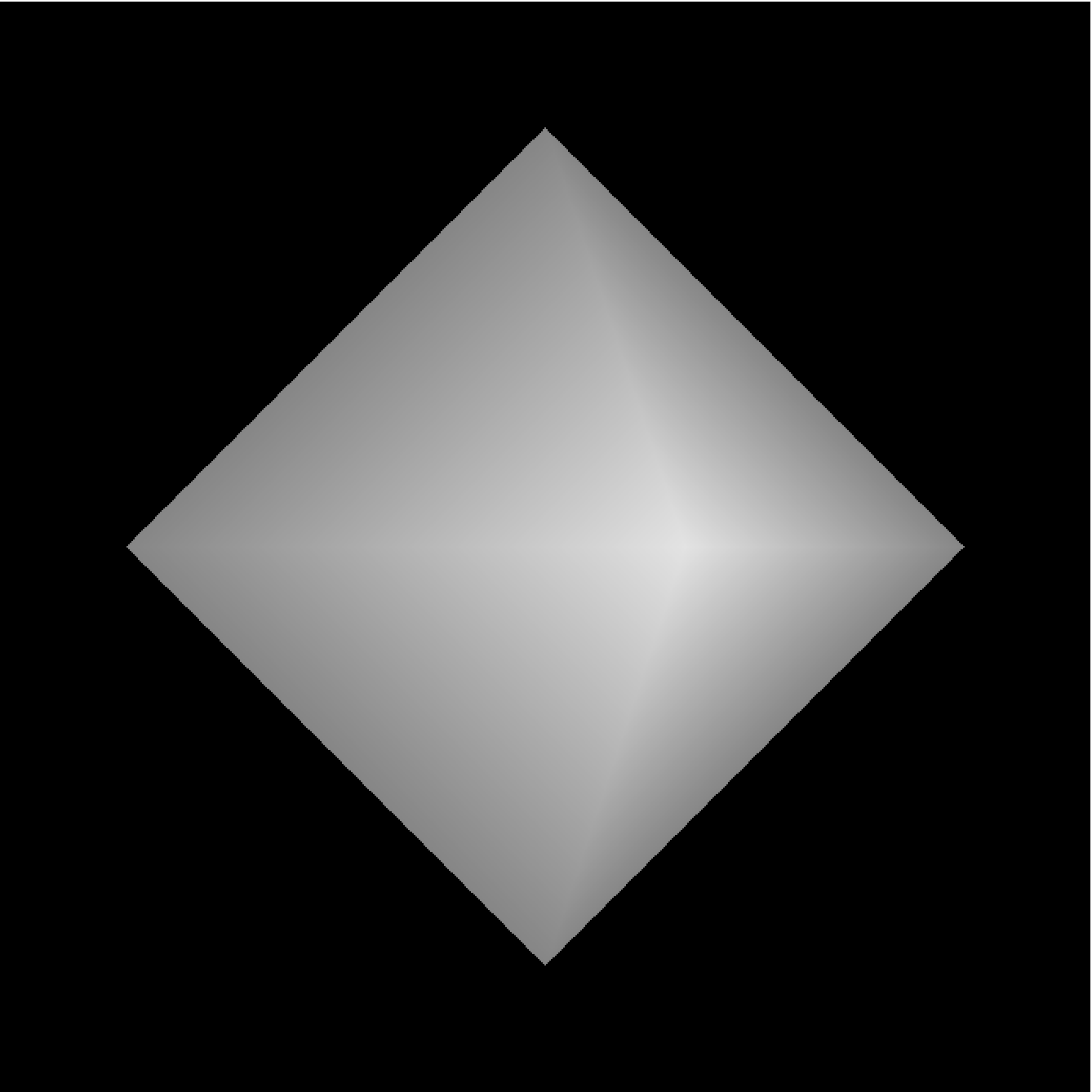}\label{fig5.8}}
\caption{Four 3D slices of the \textit{Mandelbric}}\label{figAllG2}
\end{figure}

Garant-Pelletier and Rochon \cite{GarantRochon} showed that $\mMan_3^2$ has eight principal 3D slices. So, according to \eqref{eq5.2.1} and the eight principal slices defined in \cite{GarantRochon}, we have the next lemma that corresponds to the first case discussed in the section 2, \textit{i.e.} iterates of $Q_{p,c}^m(0)$ for $m\in \mN$ are closed in the set $\mM (\bk , \bil )$.
\begin{lemma}[Parisé \cite{Parise}]\label{SymG1}
We have the following symmetries in term of $\sim$:
\begin{enumerate}
\item $\mathcal{T}^3(1,\mathbf{i_1},\mathbf{i_2}) \sim \mTet^3(1,\bk , \bil )$ $\forall \bk , \bil \in \oa \bo , \bt , \bb , \bq \fa$;
\item $\mTet^3(1, \bjp , \bjd ) \sim \mTet^3(1, \bjp , \bjt ) \sim \mTet^3(1, \bjd , \bjt )$;
\item $\mTet^3(\bo , \bt , \bjp ) \sim \mTet^3(\bk , \bil , \bk\bil )$ for all $\bk , \bil \in \oa \bo ,\bt ,\bb ,\bq \fa$ and
\item $\mTet^3(1, \bo , \bjp ) \sim \mTet^3 (1, \bk , \bjl )$ for $\bk \in \oa \bo ,\bt ,\bb , \bq \fa$ and $\bjl \in \oa \bjp , \bjd , \bjt \fa$.
\end{enumerate}
\end{lemma}
Figures \ref{fig5.1}, \ref{fig5.2}, \ref{fig5.3} and \ref{fig5.4} illustrates one slices in the four classes of 3D slices of Lemma \ref{SymG1}.
It seems that figures \ref{fig5.1} and \ref{fig5.3} looking same where these correspond to slices $\mTet^3(1,\bo , \bt )$ and $\mTet^3(\bo ,\bt , \bjp )$. Indeed, we have the next lemma that attests this remark.
\begin{lemma}\label{SymG2}
Slices $\mTet^3(1 , \bo , \bt )$ and $\mTet^3(\bo ,\bt ,\bjp )$ have the same dynamics in the sense of the relation $\sim$.
\end{lemma}
\begin{demo}
Set the numbers $c$ and $c'$ and also the function $\ph : \mM (\bo , \bt ) \rightarrow \mM (\bo , \bt )$ as
\begin{align*}
c=c_1 + c_2\bo + c_3\bt \text{,} & \quad c'=c_2\bo + c_3\bt + c_1\bjp
\end{align*}
and
\begin{equation*}
\eta =\ph (x_1 + x_2\bo + x_3 \bt + x_4 \bjp )= x_4 + x_2 \bo + x_3\bt + x_1\bjp\text{.}
\end{equation*}
So,
\begin{align*}
(\ph \circ Q_{3,c} \circ \ph^{-1} ) (\eta )&=\ph \left( (x_1^3-3x_1x_2^2-3x_1x_3^2+3x_1x_4^2 + 6x_2x_3x_4+c_1)\right.\\
& +(-x_2^3 +3x_1^2x_2 - 3x_2x_3^2 +3x_2x_4^2 - 6x_1x_3x_4 +c_2) \bo \\
& + (-x_3^3 +3x_1^2x_3 - 3x_2^2x_3 +3x_3x_4^2 - 6x_1x_2x_4 + c_3) \bt \\
&\left. + (x_4^3 + 3x_1^2x_4 - 3x_2^2x_4 - 3x_3^2x_4 + 6x_1x_2x_3)\bjp \right) \\
&=(x_4^3 + 3x_1^2x_4 - 3x_2^2x_4 - 3x_3^2x_4 + 6x_1x_2x_3)\\
&+ (-x_2^3 +3x_1^2x_2 - 3x_2x_3^2 +3x_2x_4^2 - 6x_1x_3x_4 +c_2)\bo \\
&+ (-x_3^3 +3x_1^2x_3 - 3x_2^2x_3 +3x_3x_4^2 - 6x_1x_2x_4 + c_3) \bt\\
&+(x_1^3-3x_1x_2^2-3x_1x_3^2+3x_1x_4^2 + 6x_2x_3x_4+c_1)\bjp \\
&=Q_{3,c'} \left( \eta \right)\text{.}
\end{align*}
Thus, by Definition \ref{def5.2.1}, we have the result.
\end{demo}
Because $\sim$ is an equivalence relation, by Lemmas \ref{SymG1} and \ref{SymG2}, we have find the first principal slice of $\mMan_3^3$, we will call this slice the \textit{Tetrabric}. Now, for slices that correspond to the second case (where the iterates of $Q_{3,c}^m(0)$ are not closed in $\mM (\bk , \bil )$) we have a lemma similar to Lemma \ref{SymG1}. However, when $p=3$, the iterates of $Q_{3,c}^m(0)$ are closed in $\mM (\bk ,\bil , \bim )$ (see section 2).
\begin{lemma}\label{SymG3}
We have the following symmetries:
\begin{enumerate}
\item $\mTet^3(\bo ,\bt ,\bb ) \sim \mTet^3(\bk , \bil ,\bim )$ for $\bk , \bil , \bim \in \oa \bo ,\bt , \bb , \bq \fa$;
\item Every slices of the form $\mTet^3(\bk , \bil , \bjm )$ where $\bk\bil \neq \bjm$ ,$\bk , \bil \in \oa \bo , \bt , \bb , \bq \fa$, $\bk \neq \bil$ and $\bjm \in \oa \bjp , \bjd , \bjt \fa$. Precisely,
\begin{align*}
\mathcal{T}^3(\mathbf{i_1},\mathbf{i_2},\mathbf{j_2})& \sim \mathcal{T}^3(\mathbf{i_1},\mathbf{i_2},\mathbf{j_3}) \sim \mathcal{T}^3(\mathbf{i_1},\mathbf{i_3},\mathbf{j_1}) \sim \mathcal{T}^3(\mathbf{i_1},\mathbf{i_3},\mathbf{j_3}) \sim \mathcal{T}^3(\mathbf{i_1},\mathbf{i_4},\mathbf{j_1})\\
& \sim \mathcal{T}^3(\mathbf{i_1},\mathbf{i_4},\mathbf{j_2}) \sim \mathcal{T}^3(\mathbf{i_2},\mathbf{i_3},\mathbf{j_1}) \sim \mathcal{T}^3(\mathbf{i_2},\mathbf{i_3},\mathbf{j_2}) \sim \mathcal{T}^3(\mathbf{i_2},\mathbf{i_4},\mathbf{j_1})\\
& \sim \mathcal{T}^3(\mathbf{i_2},\mathbf{i_4},\mathbf{j_3}) \sim \mathcal{T}^3(\mathbf{i_3},\mathbf{i_4},\mathbf{j_2}) \sim \mathcal{T}^3(\mathbf{i_3},\mathbf{i_4},\mathbf{j_3})\text{;}
\end{align*}
\item Every slices of the form $\mTet^3(\bk , \bjl , \bjm )$ where $\bk \in \oa \bo , \bt , \bb , \bq \fa$, $\bjl, \bjm \in \oa \bjp , \bjd , \bjt \fa$ and $\bjl \neq \bjm$. Precisely,
\begin{align*}
\mathcal{T}^3(\mathbf{i_1},\mathbf{j_1},\mathbf{j_2}) & \sim \mathcal{T}^3(\mathbf{i_1},\mathbf{j_1},\mathbf{j_3}) \sim \mathcal{T}^3(\mathbf{i_1},\mathbf{j_2},\mathbf{j_3}) \sim \mathcal{T}^3(\mathbf{i_2},\mathbf{j_1},\mathbf{j_2}) \sim \mathcal{T}^3(\mathbf{i_2},\mathbf{j_1},\mathbf{j_3}) \\
& \sim \mathcal{T}^3(\mathbf{i_1},\mathbf{j_2},\mathbf{j_3}) \sim \mathcal{T}^3(\mathbf{i_3},\mathbf{j_1},\mathbf{j_2}) \sim \mathcal{T}^3(\mathbf{i_3},\mathbf{j_1},\mathbf{j_3}) \sim \mathcal{T}^3(\mathbf{i_3},\mathbf{j_2},\mathbf{j_3}) \\
& \sim \mathcal{T}^3(\mathbf{i_4},\mathbf{j_1},\mathbf{j_2}) \sim \mathcal{T}^3(\mathbf{i_4},\mathbf{j_1},\mathbf{j_3}) \sim \mathcal{T}^3(\mathbf{i_4},\mathbf{j_2},\mathbf{j_3}) \text{ and;}
\end{align*}
\item $\mTet^3(\bjp , \bjd , \bjt )$ with itself.
\end{enumerate}
\end{lemma}
Proof of Lemma \ref{SymG3} can be found in \cite{Parise}. The same ideas from the proof of Lemma \ref{SymG1} are used in the proof of Lemma \ref{SymG3} but instead of using the set $\mM (\bk , \bil )$ we use the set $\mM (\bk ,\bil , \bim )$ to define the function $\ph$. Figures \ref{figAllG2} illustrate one slice in each four classes of 3D slices of Lemma \ref{SymG3}. From these figures, we remark that the classes of $\mTet^3(1,\bo , \bt )$ and $\mTet^3(\bo , \bt , \bjd )$ generate the same images. We notice the same phenomenon for the slices $\mTet^3(1, \bo , \bjp)$ and $\mTet^3(\bo , \bjp , \bjd )$ and also $\mTet^3(1, \bjp , \bjd )$ and $\mTet^3(\bjp , \bjd , \bjt )$. Indeed, we have this next lemma.
\begin{lemma} \label{SymG4}
We have that
\begin{enumerate}
\item $\mTet^3 (1 , \bo , \bt ) \sim \mTet^3(\bo ,\bt ,\bjd )$;
\item $\mTet^3(1 , \bo , \bjp ) \sim \mTet^3(\bo , \bjp ,\bjd )$ and;
\item $\mTet^3(1, \bjp ,\bjd ) \sim \mTet^3(\bjp , \bjd , \bjt )$.
\end{enumerate}
\end{lemma}
\begin{demo}
We prove point \textit{1)} of this lemma. Set the numbers $c$ and $c'$ as
\begin{align*}
c=c_1 + c_2\bo + c_3\bt \text{,} & \quad c'=c_2\bo + c_3\bt + c_1 \bjd \text{.}
\end{align*}
Now, let define $\ph : \mM (\bo , \bt ) \rightarrow \mM (\bo , \bt , \bjd )$ as
\begin{equation*}
\eta = \ph (x_1 +x_2\bo + x_3\bt + x_4\bjp ) = x_2 \bo + x_3\bt + x_1\bjd - x_4\bjt \text{.}
\end{equation*}
We obtain
\begin{align*}
Q_{3,c}(\ph^{-1}(\eta ) ) &= (x_1^3-3x_1x_2^2-3x_1x_3^2+3x_1x_4^2 + 6x_2x_3x_4 + c_1)\\
& +(-x_2^3 +3x_1^2x_2 - 3x_2x_3^2 +3x_2x_4^2 - 6x_1x_3x_4 +c_2) \bo \\
& + (-x_3^3 +3x_1^2x_3 - 3x_2^2x_3 +3x_3x_4^2 - 6x_1x_2x_4 + c_3) \bt \\
& + (x_4^3 + 3x_1^2x_4 - 3x_2^2x_4 - 3x_3^2x_4 + 6x_1x_2x_3)\bjp
\end{align*}
and
\begin{align*}
Q_{3,c'}(\eta )&= (-x_2^3 +3x_1^2x_2 - 3x_2x_3^2 +3x_2x_4^2 - 6x_1x_3x_4 +c_2) \bo \\
&+  (-x_3^3 +3x_1^2x_3 - 3x_2^2x_3 +3x_3x_4^2 - 6x_1x_2x_4 + c_3) \bt \\
&+ (x_1^3-3x_1x_2^2-3x_1x_3^2+3x_1x_4^2 + 6x_2x_3x_4 + c_1) \bjd \\
&- (x_4^3 + 3x_1^2x_4 - 3x_2^2x_4 - 3x_3^2x_4 + 6x_1x_2x_3) \bjt \text{.}
\end{align*}
Hence, by applying $\ph$ on the expression of $Q_{3,c}$, we have that $(\ph \circ Q_{3,c} \circ \ph^{-1}) (\eta ) = Q_{3,c'}(\eta )$ for every $\eta \in \mM (\bo , \bt , \bjd )$. Thus, $\mTet^3(1 , \bo , \bt ) \sim \mTet^3(\bo ,\bt ,\bjd )$. For the second part, set the numbers $c$ and $c'$, and also the function $\ph : \mM (\bo , \bjp ) \rightarrow \mM (\bo ,\bjp , \bjd )$ as
\begin{align*}
c=c_1+c_2\bo +c_3\bjp \text{,} & \quad c'=c_2\bo + c_3\bjp + c_1\bjd
\end{align*}
and
\begin{equation*}
\eta = \ph (x_1 + x_2\bo - x_3 \bt + x_4 \bjp ) = x_2\bo + x_4\bjp + x_1\bjd - x_3 \bq \text{.}
\end{equation*}
One can verify that $(\ph \circ Q_{3,c} \circ \ph^{-1})(\eta ) = Q_{3,c'}(\eta )$ for all $\eta \in \mM (\bo , \bjp , \bjd )$. Finally, for the last part, set the numbers $c$, $c'$ and the function $\ph : \mM (\bjp , \bjd ) \rightarrow \mM (\bjp , \bjd , \bjt )$ as
\begin{align*}
c=c_1 + c_2\bjp + c_3 \bjd \text{,} & \quad c'=c_1\bjp + c_2\bjd + c_3\bjt
\end{align*}
and
\begin{equation*}
\eta = \ph (x_1 + x_2 \bjp + x_3 \bjd - x_4 \bjt ) =-x_4 + x_1\bjp + x_2 \bjd + x_3\bjt \text{.}
\end{equation*}
Thus, $(\ph \circ Q_{3,c} \circ \ph^{-1})(\eta ) = Q_{3,c'}(\eta )$ for every $\eta \in \mM (\bjp , \bjd , \bjt )$.
\end{demo}
From the previous lemmas, we obtain the following corollary.
\begin{corollary}\label{NbPSl}
There are four principal 3D slices of the tricomplex \textit{Mandelbric}:
\begin{enumerate}
\item $\mTet^3(1, \bo , \bt )$ called \textit{Tetrabric};
\item $\mTet^3(1 , \bjp , \bjd )$ called \textit{Perplexbric};
\item $\mTet^3(1 , \bo , \bjp )$ called \textit{Hourglassbric};
\item $\mTet^3(\bo , \bt , \bb )$ called \textit{Metabric}.
\end{enumerate}
\end{corollary}
We now treat the second case of Corollary \ref{NbPSl} and we show that the \textit{Perplexbric} is an octahedron of edges $\frac{2\sqrt{2}}{3\sqrt{3}}$.

\subsection{Special case: slice $\mathcal{T}^3(1,\bjp , \bjd )$}
We had proved in section 4 that the hyperbolic \textit{Mandebric} (called the \textit{Hyperbric}) is a square of edges $\frac{2\sqrt{2}}{3\sqrt{3}}$ (see Theorem \ref{t3.2}). Now, our interest is to generalize the \textit{Hyperbric} in three dimensions. Let adopt the same notation as in \cite{GarantRochon} for the \textit{Perplexbric}
\begin{equation}
\mathcal{P}^3:=\mathcal{T}^3(1,\mathbf{j_1},\mathbf{j_2})
=\left\lbrace c=c_1+c_4\mathbf{j_1}+c_6\mathbf{j_2} \, | \, c_i \in \mathbb{R} \text{ and } \left\lbrace Q_{3,c}^m(0)\right\rbrace_{m=1}^{\infty} \text{ is bounded} \right\rbrace \text{.}\label{eq3.3.1}
\end{equation}
Before proving this result, we need this next lemma.

\begin{lemma}\label{lem3.3.1}
We have the following characterization of the \textit{Perplexbric}
\begin{equation*}
\mathcal{P}^3=\bigcup_{y\in \left[\frac{-2}{\sqrt{27}},\frac{2}{3\sqrt{3}}\right]} \left\lbrace \left[ (\mathcal{H}^3-y\mathbf{j_1})\cap (\mathcal{H}^3+y\mathbf{j_1})\right] + y\mathbf{j_2}\right\rbrace
\end{equation*}
where $\mathcal{H}^3$ is the \textit{Hyperbric} (see section 4).
\end{lemma}

\begin{demo}
By Definition of $\mathcal{P}^3$ and the idempotent representation, we have that
\begin{equation}
\mathcal{P}^3=\left\lbrace c=\left( d-c_6\mathbf{j_1}\right) \gamma_2 + \left( d + c_6\mathbf{j_1}\right) \overline{\gamma}_2 \, | \, \left\lbrace Q_{3,c}^m(0)\right\rbrace_{m=1}^{\infty} \text{ is bounded} \right\rbrace \label{PDefRemod}
\end{equation}
where $d=c_1+c_4\bjp \in \mH$. Furthermore, the sequence $\left\lbrace Q_{3,c}^m(0)\right\rbrace_{m=1}^{\infty}$ is bounded iff the two sequences $\left\lbrace Q_{3,d-c_6\mathbf{j_1}}^m(0)\right\rbrace_{m=1}^{\infty}$ and $\left\lbrace Q_{3,d+c_6\mathbf{j_1}}^m(0)\right\rbrace_{m=1}^{\infty}$ are bounded. To continue, we make the following remark about hyperbolic dynamics: $\forall z \in \mathbb{D}$
\begin{equation}
\mathcal{H}^3-z:=\left\lbrace c \in \mathbb{D} \, | \, \left\lbrace Q_{3,c+z}^m(0)\right\rbrace_{m=1}^{\infty} \text{ is bounded } \right\rbrace \text{.}\label{eq3.3.16}
\end{equation}
By Definition \ref{d3.1}, $\left\lbrace Q_{3,d-c_6\mathbf{j_1}}^m(0)\right\rbrace_{m=1}^{\infty}$ and $\left\lbrace Q_{3,d+c_6\mathbf{j_1}}^m(0)\right\rbrace_{m=1}^{\infty}$ are bounded iff $d-c_6\mathbf{j_1}, d+c_6\mathbf{j_1} \in \mathcal{H}^3$. Therefore, by \eqref{eq3.3.16}, we also have that $d-c_6\mathbf{j_1}, d+c_6\mathbf{j_1} \in \mathcal{H}^3$ iff $d \in (\mathcal{H}^3-c_6\mathbf{j_1}) \cap (\mathcal{H}^3+c_6\mathbf{j_1})$. Hence,
\begin{align*}
\mathcal{P}^3&=\left\lbrace c=c_1+c_4\mathbf{j_1}+c_6\mathbf{j_2} \, | \, c_1+c_4\mathbf{j_1} \in (\mathcal{H}^3-c_6\mathbf{j_1}) \cap (\mathcal{H}^3+c_6\mathbf{j_1}) \right\rbrace \\
&=\bigcup_{y\in \mathbb{R}} \left\lbrace \left[(\mathcal{H}^3-y\mathbf{j_1})\cap (\mathcal{H}^3+y\mathbf{j_1})\right] + y\mathbf{j_2}\right\rbrace\text{.}
\end{align*}
In fact, by Theorem \ref{t3.2},
\begin{equation}
(\mathcal{H}^3-y\mathbf{j_1})\cap (\mathcal{H}^3+y\mathbf{j_1})=\emptyset
\end{equation}
whenever $y\in \left[ -\frac{2}{3\sqrt{3}},\frac{2}{3\sqrt{3}}\right]^{c}$. This conduct us to the desire characterization of the \textit{Perplexbric}:
\begin{equation*}
\mathcal{P}^3=\bigcup_{y\in \left[\frac{-2}{\sqrt{27}},\frac{2}{3\sqrt{3}}\right]} \left\lbrace \left[(\mathcal{H}^3-y\mathbf{j_1})\cap (\mathcal{H}^3+y\mathbf{j_1})\right] + y\mathbf{j_2}\right\rbrace \text{.}
\end{equation*}
\end{demo}
As a consequence of Lemma \ref{lem3.3.1} and Theorem \ref{t3.2}, we have the following result:

\begin{theorem}\label{theo3.2.5}
$\mathcal{P}^3$ is an octahedron of edges $\frac{2\sqrt{2}}{3\sqrt{3}}$.
\end{theorem}

\section{Conclusion}
In this article, we have treated \textit{Multibrot} sets for complex, hyperbolic and tricomplex numbers. Many results presented in this article can be generalized for arbitrary integers of degree $p\geq 2$.

For the case of complex \textit{Multibrot} sets, it would be grateful if we can grade-up the proof of Theorem \ref{t2.3.1} for all Multibrots. However, as we can see, the proof is increasing in level of technicality as the degree of the polynomial $Q_{p,c}$ increases. So, we must find a different approach to prove the following next conjecture.
\begin{conjecture}\label{cj1}
Let $\mManp$ be the generalized Mandelbrot set for the polynomial $Q_{p,c}(z)=z^p+c$ where $z,c \in \mC$ and $p\geq 2$ an integer. Then, we have two cases for the intersection $\mManp \cap \mR$:
\begin{itemize}
\item[i.] If $p$ is even, then $\mManp \cap \mR = \left[-2^{\frac{1}{p-1}},(p-1)p^{\frac{-p}{p-1}} \right]$;
\item[ii.] If $p$ is odd, then $\mManp \cap \mR=\left[ -(p-1)p^{\frac{-p}{p-1}} , \, (p-1)p^{\frac{-p}{p-1}} \right]$.
\end{itemize}
\end{conjecture}
This would conduct us to another conjecture about the \textit{Hyperbrots}.
\begin{conjecture}\label{cj2}
The Hyperbrots are squares and the following characterization of \textit{Hyperbrots} holds:
\begin{enumerate}
\item If $p$ is even, then $\mHyb^p = \oa c=a+b\bj \, | \, 2^{\frac{1}{p-1}}\leq a-b, a+b \leq (p-1)p^{\frac{-p}{p-1}} \fa$;
\item If $p$ is odd, then $\mHyb^p = \oa c=a+b\bj \, | \, |a|+|b| \leq (p-1)p^{\frac{-p}{p-1}} \fa $.
\end{enumerate}
\end{conjecture}
Further explorations of 3D slices of the tricomplex \textit{Multibrot} sets are also planned. Particularly, we are interested in the case where the degree of the tricomplex polynomial is an integer $p>3$.
%-------------------------------------------------------------------------------------------------------

\section*{Acknowledgments}
DR is grateful to the Natural Sciences and
Engineering Research Council of Canada (NSERC) for financial
support. POP would also like to thank the NSERC for the award of a Summer undergraduate Research grant.

%-----------------------------------------------------------------------------------------------

\end{document}